\documentclass[a4paper]{article}
\usepackage[T1]{fontenc}
\usepackage[utf8]{inputenc}
\usepackage{framed}
\usepackage[english]{babel}
\usepackage{amssymb}
\usepackage{amsthm}
\usepackage{amsmath}
\usepackage{stmaryrd,mathtools}
\mathtoolsset{showonlyrefs}
\usepackage{bbm}
\usepackage{tikz-cd}
\usepackage{graphicx}
\usepackage[footnote,draft,danish,silent,nomargin]{fixme}
\usepackage{enumerate}
\usepackage[titletoc,title]{appendix}
\usepackage{url}

\usepackage{geometry}
\theoremstyle{plain}
\newtheorem{theorem}{Theorem}
\newtheorem{definition}{Definition}
\newtheorem{remark}{Remark}

\newtheorem{conjecture}{Conjecture}
\newtheorem{conjecture*}{Conjecture}
\newtheorem{theorem*}{Theorem}
\newtheorem*{Theorem*}{Theorem A}
\newtheorem*{Ttheorem*}{Theorem B}

\newtheorem*{ftheorem*}{Main Theorem} 
\newtheorem{lemma}{Lemma}
\newtheorem*{lemma*}{Lemma}
\newtheorem{proposition}[theorem]{Proposition}
\newtheorem*{proposition*}{Proposition}
\newtheorem{corollary}[theorem]{Corollary}
\newtheorem{corollary*}{Corollary}

\newtheorem*{ex*}{Example}

\newtheorem*{pr*}{Problem}

\newcommand{\R}{\mathbb{R}}
\newcommand{\Q}{\mathbb{Q}}
\newcommand{\C}{\mathbb{C}}

\newcommand{\N}{\mathbb{N}}
\newcommand{\Z}{\mathbb{Z}}

\newcommand\emptyarg{{}\cdot{}}

\DeclareMathOperator{\tr}{tr} 

\DeclareMathOperator{\Li}{Li} 

\DeclareMathOperator{\Spin}{Spin^c} 
\DeclareMathOperator{\im}{Im}

\DeclareMathOperator{\SU}{SU}

\DeclareMathOperator{\SL}{SL}

\DeclareMathOperator{\z}{Z}

\DeclareMathOperator{\re}{Re}

\DeclareMathOperator{\CS}{CS}

\DeclareMathOperator{\s}{S_{CS}}
\DeclareMathOperator{\dv}{{d}\hskip-.6mm}

\DeclareMathOperator{\Res}{Res}

\DeclareMathOperator{\modu}{mod}

\DeclareMathOperator{\h}{H}

\DeclareMathOperator{\I}{I}

\DeclareMathOperator{\Disc}{D}
\DeclareMathOperator{\Sph}{S}

\begin{document}

\title{Resurgence Analysis of Quantum Invariants of Seifert Fibered Homology Spheres}

\author{Jørgen Ellegaard Andersen and William Elbæk Mistegård   
}


\maketitle 

\markboth{ Jørgen Ellegaard Andersen and William Elbæk Petersen}{Resurgence Analysis of Quantum Invariants}

%
%
%
%

\begin{abstract}
	For a Seifert fibered homology sphere $X$ we show that the $q$-series invariant  $\hat{\z}_0(X;q)$ introduced by Gukov-Pei-Putrov-Vafa, is a  resummation of the Ohtsuki series $\z_0(X)$. We show that for every even $k \in \N$ there exists a full asymptotic expansion of $ \hat{\z}_0(X;q) $ for $q$ tending to $e^{2\pi i/k}$, and in particular that the limit $\hat{\z}_0(X;e^{2\pi i/k})$ exists and is equal to the  WRT  quantum invariant $\tau_k(X)$. We show that the poles of the Borel transform of $\z_0(X)$ coincide with the classical complex Chern-Simons values, which we further show classifies the corresponding  components of the moduli space of flat $\rm{SL}(2,\mathbb{C})$-connections. 
	\end{abstract}

\pagestyle{myheadings}

\date{}

\maketitle

\section*{Introduction} \label{introduction}

Let $X$ be a closed and oriented $3$-manifold and consider the level $k$ Witten-Reshetikhin-Turaev quantum invariant \cite{ReshetikhinTuraev90,ReshetikhinTuraev91}  \begin{equation} \label{eq:WRTinvariant} \tau_k(X)  \in \C, k \in \mathbb{N}.\end{equation} 
In \cite{GukovPeiPutrovVada18} Gukov, Pei, Putrov and Vafa proposed the existence of an invariant of $(X,a)$, $a \in \Spin(X),$  which is an integer power series  convergent inside the unit disc \begin{equation} \label{eq:GPPVinvariant} \hat{\z}_a(X;q) \in 2^{-c} q^{\Delta} \Z[[q]]. \end{equation}   In \cite{GukovMarinoPutrov}, $\hat{\z}_a(X;q)$ was conceived as a resummation of the asymptotic expansion of $\tau_k(X)$. The radial limit conjecture \cite{ChengChunFerrariGukovHarrison18, GukovCiprian19,GukovPeiPutrovVada18} postulates that the following limits exists $$ \hat{\z}_a(X,e^{2\pi i/k})=\lim_{q \uparrow e^{2\pi i/k}} \hat{\z}_a(X,q)$$
 and that these limits recovers $\tau_k(X)$  through an $S$-transform.

We now summarize the main results of this paper. Let $p_1,...,p_n \in \N$  be pairwise coprime integers, and let for the rest of this paper $X$ denote the oriented Seifert fibered homology sphere with $n\geq 3$ exceptional fibers
 \begin{equation} X= \Sigma(p_1,...,p_n).
 \end{equation} Let $\mathcal{M}(G)=\mathcal{M}(X,G)$ be the moduli space of flat $G$-connections, let $\s$ be the Chern-Simons action and let  $\z_0= \z_0(X)$ be the Ohtsuki series. Denote the set of classical complex Chern-Simons values by
$$\CS_{\C}(X)=\s(\mathcal{M}(\SL(2,\C))).$$ 
  \begin{itemize} 
  	  \item Theorem \ref{thm:cs} computes $\CS_{\C}(X)$, establishes that $\s$ is injective on $\pi_0(\mathcal{M}(\SL(2,\C)))$ and  that
  	$$ \pi_0(\mathcal{M}(\SL(2,\R))\sqcup_{\mathcal{M}(\text{U}(1))} \mathcal{M}(\SU(2))) \cong \pi_0(\mathcal{M}(\SL(2,\C))).$$
  \item Theorem \ref{2.1} establishes that $\CS_{\C}(X)$  coincide with the poles of the Borel transform of  $\z_0(X)$. 
  \item Theorem \ref{thm:zedhat} establishes that $\hat{\z}_0(X;q)$ is a Borel-Laplace resummation of $\z_0(X)$.
  \item  Theorem \ref{thm:radiallimit} establishes that $\hat{\z}_0(X;q)$ admits a full asymptotic expansion for $q \rightarrow e^{2\pi i/k}$, and that the radial limit conjecture is true for $X$.


  \end{itemize}
We now present these results in full detail.

\subsection*{Complex Chern-Simons Theory on $X$}

For a rational number $x \in \Q$ let $[x]= x \mod \Z.$ We prove

\begin{theorem} \label{thm:cs}
	The Chern-Simons action $\s$  is injective on $\pi_0(\mathcal{M}(\SL(2,\C)))$ and we have
	\begin{equation} \CS_{\C}(X)= \{[0]\} \sqcup\left\{ \left[ \frac{-m^2}{4P}\right]  : \text{  $m \in \Z$ is divisible by at most $n-3$ of the $p_j$'s}\right\}.
	\end{equation}
	The natural inclusion $ \mathcal{M}(\SL(2,\R))\sqcup_{\mathcal{M}(\text{U}(1))} \mathcal{M}(\SU(2)) \rightarrow  \mathcal{M}(\SL(2,\C))$
	induces an isomorphism on the level of $\pi_0$
	$$ \pi_0(\mathcal{M}(\SL(2,\R))\sqcup_{\mathcal{M}(\text{U}(1))} \mathcal{M}(\SU(2))) \cong \pi_0(\mathcal{M}(\SL(2,\C))).$$
	\end{theorem} 

\subsection*{The Borel transform and Complex Chern-Simons} Set $P=\prod_{j=1}^n p_j.$ For $k\in \N^*$ set $q_k=\exp(2\pi i/k)$ and let $\varsigma \in \C^*$ and $ \phi \in \Q$ be the constants introduced below in \eqref{eq:constants}. We consider the normalized quantum invariant
\begin{equation}
\label{def:normquant}
\widetilde{\z} _k(X)= \varsigma q_k^{\frac{\phi}{4}}  \frac{\tau_k(X)}{\tau_k(\Sph^2\times \Sph^1)} .
\end{equation} 
Let $\CS^*(X)=\CS(X)\setminus \{[0]\}$. Introduce the rational function
\begin{equation} \label{eq:G} G(z)=\prod_{j=1}^n\left(z^{\frac{P}{p_j}}-z^{-\frac{P}{p_j}}\right)\left(z^P-z^{-P}\right)^{2-n}.
\end{equation} 
In Theorem \ref{2.1} we use the notion of a resurgent function and the Borel transform, which are recalled below in Definitions \ref{def:res} and \ref{def:Borel} respectively. Let $\kappa=\sqrt{2\pi i P}.$ Building on the work of Lawrence and Rozansky\footnote{A comparison is given at the end of the introduction.}  \cite{LawrenceRozansky} we prove the following 
\begin{theorem} \label{2.1}
	There are uniquely determined polynomials $\z_{\theta} = \z_{\theta}(X)$ for $\theta \in\CS^*_{\C}(X)\}$ of degree at most $n-3$ and a formal power series $\z_{0}=\z_{0}(X) \in x^{-\frac{1}{2}}\C[[x^{-1}]]$ giving the full asymptotic expansion in the Poincare sense
	\begin{equation} \label{EXPANSIOON}
	\widetilde{\z} _k(X) \sim_{k \rightarrow \infty} \sum_{\theta \in \CS_{\C}(X)} e^{2\pi ik \theta}\z_{\theta}(k).  
	\end{equation}
	The series $\z_0$ is a normalization of the Ohtsuki series of $X$ (see Equation \eqref{eq:Ohtsuki}) whose Borel transform 
	$\mathcal{B}(\z_{0})$ is the resurgent function 
	\begin{equation} \label{ComputingBorel}
	\mathcal{B}(\z_{0})(\zeta)  = \frac{\kappa}{\pi i \sqrt{\zeta}} G\left( \exp\left(\frac{\kappa\sqrt{\zeta}}{P} \right)\right)= \frac{\kappa i}{4\pi } \frac{\prod_{j=1}^n \sinh\left(\frac{\kappa\sqrt{\zeta}}{p_j}\right)}{ \sqrt{\zeta} \left(\sinh\left(\kappa\sqrt{\zeta}\right)\right)^{n-2}},
	\end{equation} and if $\Omega$ is the set of poles of $\mathcal{B}(\z_{0})$, then
	\begin{equation} \label{Inclu}
\Omega =-2\pi i \CS_{\C}^*(X)+ 2\pi i \Z.
	\end{equation} 

\end{theorem}

\begin{remark} In accordance with the asymptotic expansion conjecture (formulated by the first author in \cite{Andersen12Crelle}, and proven in our joint work \cite{AP2} for generic mapping tori and by the first author et al. in \cite{AJHM} for finite order mapping tori with links s and further proven by Marche and Charles for certain surgeries on the figure 8 knot in \cite{CharlesMarchII}), we expect that the sum in (\ref{EXPANSIOON}) should only range over the Chern-Simons values of flat $\SU(2)$-connections 
$$ \CS(X) = \s(\mathcal{M}(X,\SU(2))),$$
e.g. that the terms which does not correspond to such values in the sum in (\ref{EXPANSIOON}) vanishes.  
This is known to be true for $n=3$  \cite{Hikami05b} and in some cases for $n=4$ \cite{Hikami05}.
However we see from (\ref{Inclu}) that the quantum invariants via resurgence determines all the Chern-Simons values of flat $\SL(2,\C)$-connections.
\end{remark}

\subsection*{A resurgence formula for the GPPV invariant $\mathbf{\hat{\z}}$ }

Let $\Delta \in \mathbb{Q}$ be given by Equation \eqref{def:delta} and consider the GPPV invariant
\begin{equation}
\hat{\z}_0(X;q) \in q^{-\Delta} \Z[[q]].
\end{equation} We recall the definition of $\hat{\z}_0(X;q)$ from \cite{GukovPeiPutrovVada18,GukovCiprian19} and recall that it is convergent for $\lvert q \rvert <1$. 

Set $m_0=P(n-2-\sum_{j=1}^np_j^{-1} ) \in \Z.$ There exists a sequence of integers $\{\chi_m\}_{m=m_0}^{\infty}$ such that for all  $z\in \C$ with $\lvert z \rvert <1$ 
\begin{equation} \label{eq:Gseries}
G(z)=(-1)^n\sum_{m=m_0}^{\infty} \chi_m z^m \in \Z[[z]].
\end{equation} Let $\mathfrak{h}$ denote the upper half-plane. Let $\tau \in \mathfrak{h}$ and set $q=\exp(2\pi i \tau)$ and  $\lambda=(-1)^n \frac{i}{2}(2P)^{-1/2}.$ Let $\varGamma=\varGamma_{+}+\varGamma_{-}$ be the oriented  unbounded contour depicted in Figure \ref{fig:ContourGamma}. We show that the GPPV invariant $\hat{\z}_0(X;q)$ is a resummation (see Appendix \ref{sec:Borel}) of the Ohtsuki series $\z_0(X)$.
\begin{theorem} \label{thm:zedhat}  
	\begin{equation} \label{eq:zedhat}
	q^{\Delta} \hat{\z}_0(X;q)= \frac{\lambda}{\sqrt{\tau}} \int_{\varGamma} \exp(-\zeta / \tau )  \mathcal{B}(\z_0(X)) (\zeta) \dv \zeta =  \sum_{m=m_0}^{\infty} \chi_m q^{\frac{m^2}{4P}}.
	\end{equation}
\end{theorem}

\begin{figure}[htp]
	\centering
	\includegraphics[scale=0.9]{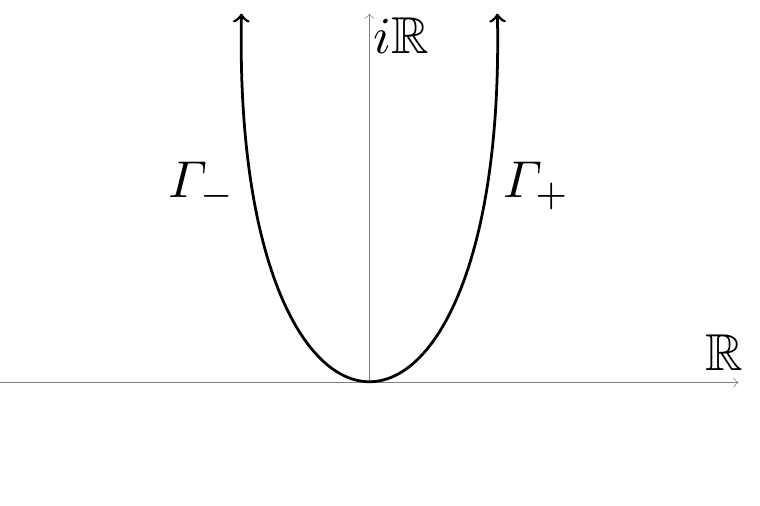}
	\caption{The integration contour $\varGamma=\varGamma_++\varGamma_{-}.$}
	\label{fig:ContourGamma}
\end{figure}  
\subsection*{The asymptotic expansion of the GPPV invariant $\hat{\z}$}
We now turn to the asymptotic expansion of $\hat{\z}_0(X)$ and the radial limit conjecture \cite{ChengChunFerrariGukovHarrison18, GukovCiprian19,GukovPeiPutrovVada18}, which is recalled  in Conjecture \ref{Conj:radiallimt}. Assume that $P$ is even. Set $\delta=\Delta-\phi/4$ and $\mu=(2\lambda \varsigma)^{-1}$. For a positive parameter $t$ set
\begin{equation}
q_{k,t}= \exp\left( \frac{2\pi i  }{k- i\frac{2Pt}{\pi }}\right) \in \mathfrak{h}.
\end{equation} 

\begin{theorem} \label{thm:radiallimit} For each $\theta \in \CS^*_{\C}(X)$  there exists a unique polynomial (defined in \eqref{def:Rtheta}) in $k$ of degree at most $n-3$ with coefficients in formal power series without constant terms
	\begin{equation}
	\check{\z}_{\theta}(k,t) \in t\cdot \Q[\pi i, k][[t]]
	\end{equation} giving a full Poincare asymptotic expansion for small $t$ and fixed even $k$
	\begin{equation} \label{eq:ZedhatExpansion}
	(\sqrt{k}\lambda)^{-1} q_{k,t}^{\Delta} \hat{\z}_0(X;q_{k,t}) \ \underset{t \rightarrow 0}{\sim} \ \tilde{\z}_k(X)+ \sum_{\theta \in \CS^*_{\C}(X)} e^{2\pi i k\theta } \check{\z}_{\theta}(k,t).
	\end{equation} In particular, for every even $k$ we have 
	\begin{equation} \label{eq:radiallimit}
	\frac{\mu }{\sqrt{k}}\lim_{q \uparrow q_{k,0}} 	q^{\delta} \hat{\z}_0(X;q) = \frac{\tau_k(X)}{\tau_k(\Sph^2 \times \Sph^1)}.
	\end{equation}
	 Thus the radial limit conjecture (Conjecture \ref{Conj:radiallimt}) holds for $X$.
	
\end{theorem}

  \subsubsection*{Comparisons with the literature}


 The existence of an asymptotic expansion
 \begin{equation}  \label{eq:LawRon} \widetilde{\z} _k(X) \sim_{k \rightarrow \infty} \sum_{ \theta \in R(Y) } e^{2\pi ik \theta} \z_{\theta}(k) \end{equation} where $R(Y) \subset \Q / \Z $
 is a finite set was proven in \cite{LawrenceRozansky}. In this work, it was also shown that $\z_0$ is a normalization of the Ohtsuki series. Our contribution in regard to \eqref{eq:LawRon} is to compute $\CS_{\C}(X)$ and to show $R(X) \subset \CS_{\C}(X).$ 
 In \cite{LawrenceRozansky} the authors do not adress the Borel transform of $\z_0$.

 The $q$-series from Theorem \ref{thm:zedhat}
 \begin{equation} \label{eq:q-series}
 \Psi(q)= \sum_{m=m_0}^{\infty} \chi_m q^{\frac{m^2}{4P}} \in \Z[[q^\frac{1}{4P}]]
 \end{equation}
 was considered in the study of  $\tau_k(\Sigma(p_1,p_2,p_3))$ by Lawrence and Zagier \cite{LawrenceZagier}, and further explored by Hikami  \cite{Hikami05}. For $n\geq 4$  Hikami in \cite{Hikami05b} considers a differently defined $q$-series. Our Theorem \ref{thm:radiallimit} generalize the result from \cite{LawrenceZagier}.
 
 The work \cite{GukovMarinoPutrov} of Gukov, Marino and Putrov is one of the main inspirations for this paper. In \cite{GukovMarinoPutrov} the authors analyse $\tau_k(X)$ for some examples with $n=3.$ The identity \eqref{Inclu} was verified for these examples. For $\tau \in \mathfrak{h}$ set $h=2\pi i\tau$ so that $q=\exp(h)$. Consider again the countor integral \begin{equation} \label{eq:contourintegral}
 \I(h)= \frac{\lambda }{ \sqrt{ \tau}} \int_{\varGamma}   \exp\left(- \frac{2\pi i \xi}{h} \right) \mathcal{B}(\z_0)(\xi) \ \dv \xi .
 \end{equation}  In \cite{GukovMarinoPutrov} identities of the form 
 \begin{equation} \label{eq:resformulaae} \I(h)=\Psi(q)
 \end{equation} where discovered. In a sense, the series $\Psi(q)$ was taken as a definition for $\hat{\z}_0(q)$ for $\Sigma(p_1,p_2,p_3)$, and the GPPV-formula (Definition \ref{df:GPPV}) was only later introduced in \cite{GukovPeiPutrovVada18}. 
 
 In the work \cite{Hiroyuki} of Fuji, Iwaki, Murakami and Terashima the $q$-series $\Psi(q)$ is also considered for general $n\geq 3$, and they prove a radial limit theorem, which is analogous to \eqref{eq:radiallimit}. The also prove an identity of the form \eqref{eq:resformulaae}. In \cite{Hiroyuki} they do not however work with the definition of the GPPV invariant $\hat{\z}_{0}(q)$, although they conjecture that this is equal to $\Psi(q)$. They also consider the  case of the WRT invariant of a knot inside $X$ and prove a difference equation for $\Psi(q)$.
 
 Our Theorem \ref{thm:zedhat} shows 
 $$q^{\Delta}\hat{\z}_0(q)=\I(h)=\Psi(q)$$ for all Seifert fibered integral homology $3$-spheres $X$ with $n\geq 3$ singular fibers, where $\hat{\z}_0(q)$ is independently defined via the GPPV-formula. We remark that those of our results that overlap with \cite{Hiroyuki} had been presented prior to their submission by the first name author in the online seminar \cite{AndersenRes2020} and by the second author at a seminar \cite{MisteRes2020} at IST, Austria. The $q$-series $\Psi(q)$ was also conjectured to be a normalization of $\hat{\z}_0(q)$ in the second author's thesis \cite{williamelbaekmistegardQuantumInvariantsChernSimons2019} . We also remark that our proof of the radial limit formula \eqref{eq:radiallimit} differs from theirs; our stronger Theorem \ref{thm:radiallimit} is derived using the resurgence formula for $\hat{\z}_0(q)$ from Theorem \ref{thm:zedhat}, whilst their proof of their radial limit theorem uses Gaussian reciprocity directly on $\Psi(q)$. We warmly thank them for cordial coordination.

 \subsubsection*{Acknowledgements} We warmly thank S. Gukov for  valuable discussions on the invariant $\hat{\z}_a(q).$

 \subsection*{Organization}
 
 In Section \ref{SeifertSection} we prove Theorem \ref{thm:cs} in several steps. Theorem \ref{thm:decompSL2C} gives a decomposition of the moduli space, Corollary \ref{cor:HUHU} computes the Chern-Simons invariants and Theorem \ref{ThmBrieskornSphere} proves that components of this moduli space are classified by their Chern-Simons value. In Section \ref{sec:BOREL} we prove Theorem \ref{2.1}. Corollary \ref{2.2} gives a resurgence formula for $\tilde{\z}_k(X)$, and Proposition \ref{thm:ExactGeneratingFunction} gives an exact formula for its generating function, verifying a special case of a conjecture of Garoufalidis \cite{Garoufalidis08}. In Section \ref{sec:Zhat} we prove Theorem \ref{thm:zedhat} and  in Section \ref{sec:radiallimit} we prove Theorem \ref{thm:radiallimit}. In Appendix \ref{sec:Borel} we present the definition of a resurgent function and the definition of the Borel transform, as well as some generalities on the relation between the Borel transform and the Laplace transform used in this paper. 
 
 \section{Complex Chern-Simons theory on $X$}  \label{SeifertSection}

Let $X$ be the oriented Seifert fibered homology $3$-sphere from the introduction. Choose $q_1,...,q_n \in \Z$ such that $(p_j,q_j)=1$ and
\begin{equation}
\sum_{j=1}^n \frac{q_j}{p_j}= \frac{1}{P}.
\end{equation}
Then $X$ has a surgery diagram as depicted in Figure \ref{fig:surgerylink} below.
\begin{figure}[htp]
	\centering
	\includegraphics{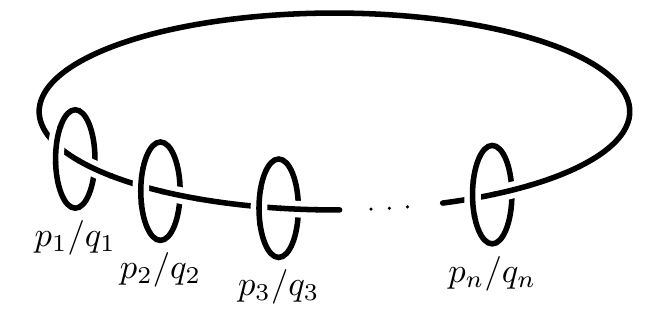}
	\caption{Surgery link for $X.$}
	\label{fig:surgerylink}
\end{figure} Without loss of generality we can assume that $p_2,...,p_n$ are odd. 
 The homeomorphism type of $X$ is unaltered under a transformation $q_j \mapsto q_j+y_j$ for any choice of integers $y_1,...,y_n$ such that $(p_j,q_j+y_j)=1$ and \begin{equation} \label{somomsum}\sum_{j=1}^n \frac{q_j}{p_j}=\sum_{j=1} \frac{q_j+y_j}{p_j}. \end{equation}
If $q_j$ is odd for $j>1$, we perform the transformation $q_j\mapsto q_j+p_j $ and $q_1 \mapsto q_1-(n-1)p_1$ which does not change the sum \eqref{somomsum}. Hence we can assume without loss of generality that $q_2,...,q_n$ are all even. 
 Recall that under our assumptions $1=P\sum_{j=1}\frac{q_j}{pj}.$ Note that this implies that $q_1$ is odd.

We recall the computation of $\tau_k(X)$ from \cite{LawrenceRozansky}. Let $S(\emptyarg,\emptyarg)$ be the Dedekind sum. Introduce the constants \begin{align} \label{eq:constants} & \varsigma  =-\frac{\sqrt{P}}{4} \exp\left(-\frac{\pi i3 }{4}\right), & \phi =3-\frac{1}{P}+12\sum_{j=1}^n S\left(\frac{P}{p_j},p_j\right). \end{align} The quantity $\phi$ is related to the Casson-Walker invariant $\lambda(X)$ \cite{Walker92} (in Casson's normalization) as follows
$$-24 \lambda(X)=\phi+P\left( n-2-\sum_{j=1}^n p_j^{-2}\right).
$$ Define the meromorphic function $F\in \mathcal{M}(\C)$ and $g \in \C[y]$  explicitly as  follows
\begin{align}  \label{def:F} \begin{split} F(y&)=\frac{1}{4}\left(\sinh\left(\frac{y}{2}\right)\right)^{2-n} \prod_{j=1}^n \sinh\left(\frac{y}{2p_j}\right), 
\\  g(y)&=\frac{iy^2}{8\pi P}. \end{split}
\end{align}  
     Lawrence and Rozansky shows in \cite{LawrenceRozansky} the following results. There exists a finite subset $R^*(X) \subset \Q^* / \Z$ and non-vanishing polynomials $\z_{\theta}(z) \in \C[z]$, $\theta \in R^*(X)$ of degree at most $n-3$ such that
\begin{equation} \label{Polypoly}
\sum_{\theta \in R^*(X)} e^{2 \pi i k \theta}\z_{\theta}(k) = - \sum_{m=1}^{2 P-1} \text{Res} \left( \frac{F(y)e^{kg(y)}}{1-e^{-ky}}, y=2\pi i m \right)
\end{equation}
for all non-negative integers $k$. Let $\gamma$ to be the contour from $(-1-i)\infty$  to $(1+i)\infty$.
Observe that $\gamma(H)$ is a steepest descent path for $g.$
Introduce the following notation
\begin{align} \label{eq:integralresiduedecomp} \begin{split}
\z^R(k) &= \sum_{\theta \in R^*(X)} e^{2 \pi i k \theta}\z_{\theta}(k)
\\ \z^{\I}(k) &= \frac{1}{2\pi i} \int_{\gamma} F(y)e^{kg(y)} \ \dv y. \end{split}
\end{align}
Recalling the definition of the normalized quantum invariant $\widetilde{\z}_k$ given in \eqref{def:normquant}, Lawrence and Rozansky proved that it can be decomposed into a sum of an integral part $\z^{\I}$ and a residue part $\z^R$
\begin{equation}
\label{phasedecomp}
\widetilde{\z} _k(X)= \z^{\I}(k)+\z^R(k).
\end{equation}
This is Equation $4.8$ on in Section $4.5$ in \cite{LawrenceRozansky}. We have used the same notation for $F,g$ and $\phi$ whereas the constant $B$ in their notation is equal to $\varsigma^{-1}$. 
Thus, if we define \begin{equation}
\z_{0}(x)= \frac{1}{2 \pi i} \sqrt{ P\pi i 8} \sum_{n=0}^{\infty} \frac{F^{(2n)}(0)\left(i8P\pi \right)^{n}}{(2n)!}  \frac{\Gamma\left(n+\frac{1}{2} \right)}{x^{n+\frac{1}{2}}} \in \ x^{-1/2}\C[[x^{-1}]]
\end{equation} 
and set $R(X)=R^*(X) \sqcup \{0\}$ then we have an asymptotic expansion
\begin{equation} \label{expansionX}
\widetilde \z_k(X) \sim_{k \rightarrow \infty} \sum_{\theta\in R(X)} e^{2\pi i k \theta}  \z_{\theta}(k).
\end{equation}
  In the work \cite{LawrenceRozansky} it was observed that $\z_{0}$ is in fact a normalization of the Ohtsuki-series \cite{Ohtsuki95,Ohtsuki96a,Ohtsuki96b}. Let $\z_{\infty}$ denote the Ohtsuki series (with the normalization used in \cite{LawrenceRozansky}). Introduce the variable $h=q_k-1$, where as above $q_k=\exp(2\pi i/k)$. In section $4.6$ of \cite{LawrenceZagier}  they show the following identity
  \begin{equation} \label{eq:Ohtsuki}
\z_0(k) = \z_{\infty}(h).
  \end{equation}
  
\subsection{The moduli space and complex Chern-Simons values}

We now begin our investigation of $\mathcal{M}(X,\SL(2,\C)),$ which closely follows \cite{FintushelStern90}. We have the following presentation of the fundamental group of $X$
\begin{equation} \label{presentation}
\pi_1(X) \simeq \left\langle h,  x_1,...,x_n \mid  x_1x_2\cdots x_n, x_j^{p_j}h^{-q_j}, [x_j,h], j=1,...,n  \right\rangle.
\end{equation}

Let us first recall a few of Fintushel and Stern's results concerning the moduli space $\mathcal{M}(X,\SU(2))$  establish in \cite{FintushelStern90}. As $X$ is an integral homology sphere, the only reducible representation into $SU(2)$ is the trivial one. For an irreducible representation $\rho:\pi_1(X) \rightarrow \SU(2)$ at most $n-3$ of the $\rho(x_j)$ are $\pm I,$ and if exactly $n-m$ of the $\rho(x_j)$ are equal to $\pm I,$ then the component of $\rho$ in $\mathcal{M}(X,\SU(2))$ is of dimension $2(n-m)-6.$

 Let $L(p_1,...,p_n) \subset \mathbb{N}^{n}$ be the the set of $n$-tuples $l=(l_1,...,l_n)$ which satisfies the following condition. We have $0\leq l_1 \leq p_1$ and $0\leq l_j \leq \frac{p_j-1}{2},$ for $j=2,...,n$ and there exists at least three distinct ${j_1}<{j_2}<{j_3}$ with $l_{j_t}\not=0$ for $t=1,2,3.$
The following proposition is an adaptation and generalization of Lemma $2$ in \cite{BodenCurtis06} and Lemma $2.1$ and Lemma  $2.2$ in \cite{FintushelStern90}. 

 \begin{proposition} \label{rho_l}
 Let $l=(l_1,l_2,...,l_n) \in L(p_1,...,p_n).$ Then there exists matrices $Q_j \in \SL(2,\mathbb{C})$ and a representation 
$ \rho_l: \pi_1(X) \rightarrow \SL(2,\mathbb{C})$  with
\begin{equation} \label{X_1}
\rho_l(x_1)= Q_1\begin{pmatrix} e^{\frac{\pi i l_1}{p_1}} & 0
\\ 0 &  e^{-\frac{\pi i l_1}{p_1}} \end{pmatrix}Q_1^{-1}, \ \ \   \rho_l(x_j)= Q_j\begin{pmatrix} e^{\frac{2\pi i l_j}{p_j}} & 0
\\ 0 &  e^{-\frac{2\pi i l_j}{p_j}} \end{pmatrix}Q_j^{-1}
\end{equation}
for $j=2,\ldots n$. In fact we can choose $Q_j=I$ for $j\not=j_2,j_3.$ 

Furthermore we can choose the $Q_j$'s such that 
$$ \rho_l : \pi_1(X) \rightarrow \SL(2,\mathbb{R})$$ 
or
$$ \rho_l: \pi_1(X) \rightarrow \SU(2)$$ 
depending on properties of $l$.
For any non-trivial representation $\rho: \pi_1(X) \rightarrow \SL(2,\mathbb{C})$ there exists $l' \in \mathbb{N}^n$ such that $p_j$ divides $ l'_j$ for at most $n-3$ of the indices $j=1,2,3,...,n$ and such that $\rho$ is of the form 
\begin{equation} \label{randomrep}
\rho(x_1)= S_1\begin{pmatrix} e^{\frac{\pi i l'_1}{p_1}} & 0
\\ 0 &  e^{-\frac{\pi il'_1 }{p_1}} \end{pmatrix}S_1^{-1}, \ \ \   \rho(x_j)= S_j\begin{pmatrix} e^{\frac{2\pi i l'_j}{p_j}} & 0
\\ 0 &  e^{-\frac{2\pi i l'_j}{p_j}} \end{pmatrix}S_j^{-1}
\end{equation}
for some $S_1,...,S_n \in \SL(2,\mathbb{C}).$

Finally we have that the map which associates $l\in L(p_1,...,p_n)$ to a non-trivial representation $\rho : \pi_1(X) \rightarrow \SL(2,\C)$ via \eqref{randomrep} induces an isomorphism
$$\pi_0(\mathcal{M}^*(X,\SL(2,\C))) \cong L(p_1,...,p_n)$$
\end{proposition}

The family of Brieskorn integral homology spheres ($n=3$) is very special due to the fact that the moduli space $\mathcal{M}(\Sigma(p_1,p_2,p_3), \SL(2,\mathbb{C}))$ is finite with cardinality given by the $\SL(2,\mathbb{C})$ Casson invariant introduced by Curtis \cite{Curtis01,Curtis03} $$\lambda_{\SL(2,\mathbb{C})}(\Sigma(p_1,p_2,p_3))=(p_1-1)(p_2-1)(p_3-1)/4.$$ This is shown by Boden and Curtis \cite{BodenCurtis06}. Prior to this and in relation to Floer homology, Fintuschel and Stern \cite{FintushelStern90} analyzed the $\SU(2)$ moduli space $\mathcal{M}(X,\SU(2))$ of the Seifert fibered three manifold $X$ considered in this paper and their work shows that the components are even dimensional manifolds with top dimension $2n-6.$ This is in stark contrast to the finiteness of the moduli space $\mathcal{M}(\Sigma(p_1,p_2,p_3), \SL(2,\mathbb{C}))$.  In the three fibered case Kitano and Yamaguchi \cite{KitanoYamaguchi16} has gives a decomposition
\begin{equation}
\mathcal{M}(\Sigma,\SL(2,\mathbb{C}))=\mathcal{M}(\Sigma,\SL(2,\mathbb{R})) \bigcup_{\mathcal{M}(\Sigma,\text{U}(1))} \mathcal{M}(\Sigma,\SU(2)).
\end{equation}
Where $\Sigma=\Sigma(p_1,p_2,p_3).$ Here we can observe the following generalization of this work as an imidiate corollary of Proposition \ref{rho_l}.

\begin{theorem} \label{thm:decompSL2C}
The natural inclusion
$$ \mathcal{M}(X,\SL(2,\R))\sqcup_{\mathcal{M}(X,\text{U}(1))} \mathcal{M}(X,\SU(2)) \rightarrow  \mathcal{M}(X,\SL(2,\C))$$
induces an isomorphism on the level of $\pi_0$
$$ \pi_0(\mathcal{M}(X,\SL(2,\R))\sqcup_{\mathcal{M}(X,\text{U}(1))} \mathcal{M}(X,\SU(2))) \cong \pi_0(\mathcal{M}(X,\SL(2,\C))).$$
\end{theorem}
By this corollary we can in particular conclude that all Chern-Simons values are real and they only depend on $l\in  L(p_1,...,p_n)$. In Proposition \ref{CSvalues} below we actually provide an explicit formula.

Before commencing the proof let us introduce the following notation
\begin{equation} \exp(x)=\begin{pmatrix} e^{x} & 0
\\ 0 &  e^{-x} \end{pmatrix}, \end{equation}
which should not cause any ambiguities as long as the context shows that we are dealing with a matrix.

\begin{proof} We start with the construction of $\rho_l.$ Introduce the matrices $$X_1=\exp(\pi i l_1/p_1), \ \ \ X_j=\exp(2\pi il_j/p_j)$$
for $j\in \{2,....,n\}\setminus \{j_2,j_3\}.$ 
Rewrite the relation $\prod_{j=1}^n x_j=1$ as the equivalent relation $$x_{{j_3}+1}\cdots x_n x_1\cdots x_{j_1} \cdots x_{j_2}\cdots x_{j_{3}-1}=x_{j_3}^{-1}.$$ 
 Assume we can chosen $Q_{j_2},Q_{j_3} \in \SL(2,\mathbb{C})$ such that
\begin{equation} \label{veryfinnny}
X_{{j_3}+1}\cdots X_n X_1\cdots X_{j_1} \cdots Q_{j_2}X_{j_2}Q_{j_2}^{-1}\cdots X_{j_{3}-1}=Q_{j_3}^{-1}X_{j_3}^{-1}Q_{j_3}.
\end{equation}
Taking $Q_j=I$ for $j\notin \{j_2,j_3\},$ we can define $\rho:\pi_1(X)\rightarrow \SL(2,\mathbb{C})$ by the assigment
\begin{equation}
\rho(x_j)=Q_jX_jQ_j^{-1}, \ \rho(h)=X_1^{p_1}.
\end{equation}
To see this, observe that $B:=X_1^{p_1}=(-I)^{l_1}$ is central and  as $q_1$ is odd whereas $q_j$ is even for $j\geq 2,$ we also have  $X_j^{p_j}=B^{q_j}, \forall j.$ The last relation in $\pi_1(X)$ is ensured by  \eqref{veryfinnny}. Observe that it will suffice to choose $Q_{j_2} \in \SL(2,\mathbb{C})$ such that
\begin{equation} \label{radadadadadadada}
 \tr \left( X_{{j_3}+1}\cdots X_n X_1\cdots X_{j_1} \cdots Q_{j_2}X_{j_2}Q_{j_2}^{-1}\cdots X_{j_{3}-1} \right) =2  \cos\left( \frac{2\pi l_{j_3}}{p_{j_3}}\right) 
\end{equation}
because this will ensure that there exists some $Q_{j_3} \in \SL(2,\mathbb{C})$ with the property that
$$Q_{j_3} X_{{j_3}+1}\cdots X_n X_1\cdots X_{j_1} \cdots Q_{j_2}X_{j_2}Q_{j_2}^{-1}\cdots X_{j_{3}-1} Q_{j_3}^{-1}=X_{j_3},$$ 
since non-diagonalizable elements of $ \SL(2,\mathbb{C})$ have trace $\pm2$, given that the unit determinant condition implies that the unique eigenvalue with multiplicity two for such elements must be either $1$ or $-1$ and we have that
$$2  \left|\cos\left( \frac{2\pi l_{j_3}}{p_{j_3}}\right) \right|<2.$$
For \eqref{radadadadadadada} we used our assumption on $j_3.$ Write
\begin{align}
&X_{{j_3}+1}\cdots X_n X_1\cdots X_{j_1} \cdots X_{j_2-1}=\exp(ia),
\\ &X_{j_2}=\exp(ib),
 \\ &X_{j_2+1} \cdots X_{j_{3}-1}=\exp(ic),
 \\ &X_{j_3} =\exp(id).
\end{align}
We observe that by the conditions on $(l_{j_2},p_{j_2})$ and $(l_{j_3},p_{j_3})$ we have that 
$$b,d \notin \pi \Z.$$
Let $Q_{j_2}=\tilde Q$ where
$$ \tilde Q =  \tilde Q(u,v,w,z) =\begin{pmatrix} u & v
\\ w & z
\end{pmatrix}$$ for $u,v,w,z$ to be chosen below. Assume $uz-vw=1$ so that $\tilde Q\in \SL(2,\mathbb{C}).$ We compute
\begin{multline}
X_{{j_3}+1}\cdots X_n X_1\cdots X_{j_1} \cdots X_{j_2-1} \tilde QX_{j_2} \tilde Q^{-1} X_{j_2+1}\cdots X_{j_{3}-1} 
\\ = \begin{pmatrix} zue^{i(a+b+c)}-wve^{i(a-b+c)} & -uve^{i(a+b-c)}+uve^{i(a-b-c)}
\\ zve^{i(b+c-a)}-zwe^{i(c-a-b)} & -vwe^{-i(a+b+c)}+zue^{i(b-a-c)} \end{pmatrix}.
\end{multline}
Thus we have that
\begin{multline} \text{tr}\begin{pmatrix} zue^{i(a+b+c)}-wve^{i(a-b+c)} & -uve^{i(a+b-c)}+uve^{i(a-b-c)}
\\ zwe^{i(b+c-a)}-zwe^{i(c-a-b)} & -vwe^{-i(a+b+c)}+zue^{i(b-a-c)} \end{pmatrix}\\ =2zu\cos(a+b+c)-2wv\cos(a+c-b).
\end{multline}
It follows that we must solve
\begin{equation} \label{system}
\begin{pmatrix} \cos(a+b+c) & \cos(a+c-b)
\\ 1 & 1
\end{pmatrix} \begin{pmatrix} zu \\ -wv \end{pmatrix}=\begin{pmatrix} 2\cos(d) \\ 1 \end{pmatrix}.
\end{equation}
 Using the trigonometric identity $\cos(x+y) - \cos(x-y)=-2\sin(x)\sin(y)$ we get
 \begin{align}
 \det \begin{pmatrix} \cos(a+b+c) & \cos(a+c-b)
\\ 1 & 1
\end{pmatrix}&= \cos((a+c)+b)-\cos((a+c)-b)
\\ &=-2\sin(a+c)\sin(b).
 \end{align}
 Thus it remains to argue $a+c \notin \pi \mathbb{Z}$ and $ b \notin \pi\mathbb{Z}.$ Assume towards a contradiction that $a+c =\pi m$ for some $m \in \mathbb{Z}.$ Hence we would have $P(a+c)=Pm\pi$ for some integer $m$, which would imply 
 \begin{equation}
 l_{j_1} 2^{\epsilon} \prod_{t\not=j_1} p_t = 0 \ \modu \ p_{j_1}
 \end{equation}
 for $\epsilon \in \{0,1\},$ with $\epsilon =0$ for $j_1=1$ and $\epsilon=0$ otherwise. This is a contradiction, as $2^{\epsilon} \prod_{t\not=j_1} p_t$ is invertible in $\mathbb{Z}/p_{j_1} \mathbb{Z}$ and $1\leq l_{j_1} \leq (p_{j_1}-1)/2^{\epsilon}.$ We see that $b\notin \pi \mathbb{Z}$ directly from the conditions on $l_{j_2}$. Thus we can solve \eqref{system}, and hence find the needed $u,v,w,z$, which concludes the proof of the first part of the proposition.
 
Let us now prove that we can actually choose the $Q_j$'s such that we obtain an $\SL(2,\R)$-representation. We will denote this new choice of the $Q_j$ by $ Q_j^\R$. We set $Q^\R_j = Q$ for $j\in \{1, \ldots j_2-1 ,j_3+1, \ldots, n\}$, where 
$$Q=\frac{\sqrt{2}}{2}\begin{pmatrix} 1& -i
\\ -i & 1
\end{pmatrix}.$$ 
which has the following property
$$ Q \exp(a) Q^{-1} = \begin{pmatrix} \cos(a) & -\sin(a)
\\ \sin(a) & \cos(a)
\end{pmatrix}.$$
Introduce the notation $\tilde X_j = Q X_j Q^{-1}$ and observe by the above computations that
\begin{multline}
\tr(\tilde X_{{j_3}+1}\cdots \tilde X_n \tilde X_1\cdots \tilde X_{j_1} \cdots \tilde X_{j_2-1}Q \tilde QX_{j_2}(Q\tilde Q)^{-1}) =
\\ \tr(X_{{j_3}+1}\cdots X_n X_1\cdots X_{j_1} \cdots X_{j_2-1} \tilde QX_{j_2}\tilde Q^{-1}) =
\\ =2zu\cos(a+b)-2wv\cos(a-b).
\end{multline}
To understand which values, say $t$, this trace can take, we consider in analogy with \eqref{system} the equation
\begin{equation} \label{system2}
\begin{pmatrix} \cos(a+b) & \cos(a-b)
\\ 1 & 1
\end{pmatrix} \begin{pmatrix} zu \\ -wv \end{pmatrix}=\begin{pmatrix} t \\ 1 \end{pmatrix}.
\end{equation}
The determinant is $D = -2\sin(a)\sin(b)$, which is non-vanishing since $a,b\notin \pi \Z$. Then we have that
\begin{equation} \label{system3}
\begin{pmatrix} zu \\ -wv \end{pmatrix} = \frac1D \begin{pmatrix} t - \cos(a-b) \\ -t+\cos(a+b) \end{pmatrix}.
\end{equation}
For $t\in \R$ we observe that $zu\in \R$ and $wv\in \R$. Now we compute
\begin{multline}
Q\tilde Q X_{j_2}(Q\tilde Q)^{-1}=
\\ \begin{pmatrix} (uz-wv)\cos(b) + (uv+wz) \sin(b) & 
(-(uz+wv) - i(uv-wz))\sin(b)
\\ ((uz+wv) - i(uv-wz))\sin(b) &
(uz-wv)\cos(b) - (uv+wz) \sin(b) \end{pmatrix}.
\end{multline}
From which we see that 
$$Q\tilde QX_{j_2}(Q\tilde Q)^{-1}\in SL(2,\R)$$ 
if and only if
$$ \text{Im}(uv+wz) = 0, \ \ \ \text{Re}(uv-wz) =0$$
or equivalently
$$ uv = \overline{wz}$$
whenever $t\in \R$. But then this implies that
$$ |u|^2 |v|^2 = \frac{(t-\cos(a+b))(t-\cos(a-b))}{D^2}.$$
Which we can solve when $t> \cos(a\pm b)$ by letting
$$ u= \frac{\sqrt{t-\cos(a+b)}}{D}, \ \ \ v= \frac{\sqrt{t-\cos(a-b)}}{D}$$
and then
$$ z = \sqrt{t-\cos(a+b)}, \ \ \ w= \sqrt{t-\cos(a-b)}$$
and when $t< \cos(a\pm b)$ then we can take
$$ u= \frac{\sqrt{\cos(a+b) - t }}{D}, \ \ \ v= \frac{\sqrt{\cos(a-b)-t}}{D}$$
and then
$$ z = \sqrt{\cos(a+b)-t}, \ \ \ w= \sqrt{\cos(a-b)-t}.$$
This allows us to complete the construction as follows. First we assume that  $a, c\notin \pi \Z$. For 
$$ \tilde Q_i = \tilde Q(u_i,v_i,w_i,z_i)$$
for $i=1,2$ we consider the equation
$$
\left[\tilde X_{{j_3}+1}\cdots \tilde X_n  \tilde X_1\cdots \tilde X_{j_1} \cdots Q \tilde Q_1X_{j_2}(Q \tilde Q_1)^{-1}\right]  =
\left[\left(\tilde X_{j_2+1} \ldots \tilde X_{j_3-1} Q \tilde Q_2 X_{j_3}(Q \tilde Q_2)^{-1}\right)^{-1}\right]
$$
which is equivalent to 
$$\tr\left(\tilde X_{{j_3}+1}\cdots \tilde X_n  \tilde X_1\cdots \tilde X_{j_1} \cdots Q \tilde Q_1X_{j_2}(Q \tilde Q_1)^{-1}\right) =
\tr\left(\tilde X_{j_2+1} \ldots \tilde X_{j_3-1} Q \tilde Q_2 X_{j_3}(Q \tilde Q_2)^{-1}\right)$$
since these are certainly all $\SL(2,\C)$ matrices.
But now, using that we also have that $c,d \notin \pi\Z$, we can chose $t$ bigger than $\cos(a \pm b)$ and $\cos(c\pm d)$ and fix $\tilde Q_i$ as above such that
\begin{equation}\label{tra1}
 \tr\left(\tilde X_{{j_3}+1}\cdots \tilde X_n  \tilde X_1\cdots \tilde X_{j_1} \cdots Q \tilde Q_1X_{j_2}(Q \tilde Q_1)^{-1}\right)  = t
 \end{equation}
and 
\begin{equation}\label{tra2}
 \tr\left(\tilde X_{j_2+1} \ldots \tilde X_{j_3-1} Q \tilde Q_2 X_{j_3}(Q \tilde Q_2)^{-1}\right) = t
 \end{equation}
Thus we can now conclude that there exist $Q^{\R} \in \SL(2,\R)$ such that
\begin{multline}\tilde X_{{j_3}+1}\cdots \tilde X_n  \tilde X_1\cdots \tilde X_{j_1} \cdots Q \tilde Q_1X_{j_2}(Q \tilde Q_1)^{-1} = \\Q^\R \left(\tilde X_{j_2+1} \ldots \tilde X_{j_3-1} Q \tilde Q_2 X_{j_3}(Q \tilde Q_2)^{-1}\right)^{-1} (Q^\R)^{-1} =
\\  Q^\R Q \tilde Q_2 X_{j_3}^{-1}  (Q^\R Q \tilde Q_2)^{-1} Q^\R  \tilde X_{j_3-1}^{-1} (Q^\R)^{-1} \ldots Q^\R  \tilde X_{j_2+1}^{-1} (Q^\R)^{-1}  .
\end{multline}
Thus if we further set 
$$ Q^\R_{j_2} = Q \tilde Q_1, Q^\R_{j_3} = Q^\R Q \tilde Q_2$$
and
$$ Q^\R_j = Q^\R Q$$
for $j\in \{j_2+1,\ldots,j_3-1\}$, then
we have found the needed conjugation to obtain an $\SL(2,\R)$-representation. Let us now consider the remanning cases.
Suppose that $a\in \pi Z$ but $c\notin \pi \Z$. Then the common trace $t$ is by \eqref{tra1} forced to be $e^{ia} \cos b$, so we can solve \eqref{tra2} if and only if $e^{ia} \cos b$ is not contained in the interval spanned by the two values $\cos(c\pm d)$. If this is the case we proceed with the argument as above. If on the other hand $e^{ia} \cos b$ is contained in the interval spanned by $\cos(c\pm d)$, then it is well known that we can choose $Q_j \in \SU(2)$ so as to obtain an $\SU(2)$-representation. A similar argument of course works in the case where $c\in \pi Z$ but $a\notin \pi \Z$. If we have $a,c \in \pi\Z$, then $a+c \in \pi\Z$, but this we have already argued is impossible.

 Now let $\rho: \pi_1(X) \rightarrow \SL(2,\mathbb{C})$ be an arbitrary non-trivial representation. As remarked before any non-trivial representation is irreducible since $X$ is an integral homology three sphere. Since $\rho(h)$ commutes with the image of $\rho$, we see that $\rho(h) = \pm I$. Hence the relation $x_j^{p_j}=h^{q_j}$ implies that $\rho(x_j)^{p_j}=\pm I,$ and for $j=2,...,n$ we must have $\rho(x_j)^{p_j}=I,$ since $q_j$ is even. Hence $\rho$ must be of the form \eqref{randomrep} for some $l'\in \mathbb{N}^n.$ It only remains to argue that at most $n-3$ of the $\rho(x_j)$ are $\pm I.$ If not, the relation $x_1x_2 \cdots x_n=1$ implies that there is $j_1<j_2$ with 
 $\rho(x_{j_1})\rho(x_{j_2})=\pm I.$ As $p_{j_1}$ and $p_{j_2}$ are relatively coprime, this is only possible if $\rho(x_{j_1})=\pm I$ and $\rho(x_{j_2})=\pm I.$ This would imply that $\rho(\pi_1(X))\subset \{\pm 1\}=Z(\SU(2))$ which contradicts the fact that $\rho$ is irreducible since it was assumed non-trivial.

 We describe the connected components of $\mathcal{M}^*(X,\SL(2,\C))$. First we assume that $l_1 \in \{1, \ldots p_1-1\}$ and $l_2 >0$ for an $l\in L(p_1,...,p_n) $. We will now prove that the subset $\mathcal{M}_l^*(X,\SL(2,\C))$ of $\mathcal{M}^*(X,\SL(2,\C))$ consisting of conjugacy classes of representations $\rho$ for which 
 $$ \tr(\rho(x_1)) = 2\cos(\pi l_1/p_1), \ \ \ \tr(\rho(x_j)) = 2\cos(2\pi l_j/p_j), \ \ j = 2, \ldots n$$
 is connected.
 Let 
 \begin{equation}
T  =    \bigtimes_{j=3}^{n} \left[\exp\left(\frac{2\pi il_j}{p_j}\right)\right].
\end{equation}
 It is obvious that $T$ is connected.  Let now $m : T \rightarrow \SL(2,\mathbb{C})$ be the algebraic product map. Let $P\subset \SL(2,\mathbb{C})$ be the set of non-diagonalizable elements in $\SL(2,\mathbb{C})$ and observe that $P$ has complex co-dimension one, thus so does $m^{-1}(P) \subset T$, but then it follows that $T' = T- m^{-1}(P)$ is also connected. For any
 $$(M_3, \ldots, M_n) \in T',$$
 we observe that the set of $(Q_1,Q)\in \SL(2,\C)$ which solves
 \begin{equation}\label{eqa1}
 \exp(\pi l_1/p_1) Q_1 \exp(2 \pi  l_2/p_2) Q_1^{-1} = Q M_n^{-1} \ldots M_3^{-1} Q^{-1}
 \end{equation}
 is non empty and connected since it is acted transitively on by $(\C^*)^2 \times \C^*$, where the first factor comes from the ambiguity from solving \eqref{system2} and the second comes from the stabiliser of $M_n^{-1} \ldots M_3^{-1}$ under conjugation. But then we see that an open dense subset of $\mathcal{M}_l^*(X,\SL(2,\C))$ is connected, thus $\mathcal{M}_l^*(X,\SL(2,\C))$ it self must be connected.
If $l_1\in \{0,p_1\}$ or $l_2 =0$, we proceed as follows. Choose $j_1, j_2 \in \{1,\ldots, n-2\}$ such that
$$ a = \pi \frac{l_1}{p_1} + 2\pi \sum_{j=2}^{j_1} \frac{l_j}{p_j}, \ \ \ b = 2\pi \sum_{j=j_1+1}^{j_2} \frac{l_j}{p_j},$$
has the property that $a,b\notin \pi \Z$. Now consider the equation
\begin{equation}\label{eqa2}
\exp(a) Q_1 \exp(b) Q_1^{-1} = Q M_n^{-1} \ldots M_{j_2+1}^{-1} Q^{-1}.
\end{equation}
The connectedness is now argued in exactly the same way, with \eqref{eqa2} in place of \eqref{eqa1}.
 \end{proof}

For an $\SL(2,\C)$-connection $a$ in the trivial $\SL(2,\C)$-bundle on $X$ we recall the Chern-Simons action is given by 
\begin{equation} \label{eq:CSaction} \s(a)=\frac{1}{8\pi^2 }\int_X \tr\left(a\dv a +\frac{2}{3}a^3\right) \mod \Z.
\end{equation} 

We now compute the Chern-Simons values of the representations constructed above. 
\begin{proposition} \label{CSvalues}  For any representation $\rho : \pi_1(X) \rightarrow SL(2,\C)$, define $l=(l_1,...,l_n) \in L(p_1,...,p_n),$ so that
$$ \tr(\rho(x_1)) = 2\cos(\pi l_1/p_1), \ \ \ \tr(\rho(x_j)) = 2\cos(2\pi l_j/p_j), \ \ j = 2, \ldots n.$$
Then we have that
\begin{equation} \label{CS}
\s\left( \rho \right) = -\frac{P}{4}\left(\frac{l_1}{p_1}+\sum_{j=2}^{n} \frac{2l_j}{p_j}\right)^2 \  \modu \ \mathbb{Z}.
\end{equation}
\end{proposition}
Formula \eqref{CS} was proven for $\SU(2)$ connections by Kirk and Klassen and it is stated in Theorem $5.2$ in \cite{KirkKlassen90}. It is proven using the following general result.
Let $M$ be a closed oriented three manifold containing a knot $K.$ Let $Y$ be the complement of a tubular neighborhood of $K$ in $M.$ With respect to an identification $M\setminus Y \simeq D^2 \times S^1,$ choose simple closed curves $\mu, \lambda$ on $\partial Y$ intersecting in a single point such that $\mu$ bounds a disc of the form $D^2 \times \{1\}.$ Let $\rho_t:\pi_1(Y) \rightarrow \text{SU}(2)$ be a path of representations such that $\rho_0(\mu)=\rho_1(\mu)=1,$ and for which there exists continuous piecewise differentiable functions 
$$\alpha,\beta:I \rightarrow \R$$ 
with 
\begin{align}
\begin{split}
\rho_t(\mu)= \begin{pmatrix} e^{2\pi i \alpha(t)} & 0
\\ 0 &  e^{-2\pi i \alpha(t)} \end{pmatrix}, \  \rho_t(\lambda)=\begin{pmatrix} e^{2\pi i \beta(t)} & 0
\\ 0 &  e^{-2\pi i \beta(t)} \end{pmatrix}.
\end{split}
\end{align}
Thinking of $\rho_1,\rho_0$ as flat connections on $M$ we have
\begin{equation} \label{KK}
\s(\rho_0)-\s(\rho_1)=-2 \int_0^1 \beta(t) \alpha'(t) \ \dv t \ \text{mod} \ \mathbb{Z}. 
\end{equation}
Notice that formula \eqref{KK} differs from the corresponding formula in \cite{KirkKlassen90} by a sign. This discrepancy was already discussed by Freed and Gompf in \cite{FreedGompf91} and is due to a sign convention. See the footnote on page $98$ in \cite{FreedGompf91}. The formula \eqref{KK} was also used in the work \cite{AndersenHansen} by the first author and Hansen.

\begin{proof}[Proof of Proposition \ref{CSvalues}]
Let $K\subset X$ be the $n$'th exceptional fiber. Let $Y$ be the complement of a tubular neighborhood of $K$ in $X.$ Removing $K$ has the effect on $\pi_1$ of removing the relation $x_n^{p_n}=h^{-q_n}$, i.e. we have a presentation
\begin{equation} \label{presentation2}
\pi_1(Y) \simeq \left\langle h,  x_1,...,x_n \mid  x_1x_2\cdots x_n, x_1^{p_1}h^{-q_1},...,x_{n-1}^{p_{n-1}}h^{-q_{n-1}} , \forall j [x_j,h]   \right\rangle.
\end{equation}
As the meridian and longitude of $\partial Y$ we can take $\mu=x_n^{p_n}h^{q_n}$ and $\lambda=x_n^{-p_1\cdots p_{n-1}}h^c$ respectively, where $c=\sum_{j=1}^{n-1} \frac{p_1\cdots p_{n-1}q_j}{p_j}.$ These choices of meridian and longitude coincide with the choices made in \cite{FintushelStern90}.

Let $\rho:\pi_1(X)\rightarrow \SL(2,\mathbb{C})$ be any irreducible representation with its corresponding $l=(l_1,...,l_n) \in L(p_1,...,p_n)$. Now $\rho(h)=\exp(\pi i v)$ for some integer $v\in \Z$.
 Introduce the two quantities 
 $$\xi= P\left(\frac{l_1}{p_1}+\sum_{j=2}^{n} \frac{2l_j}{p_j}\right)$$ and $$\eta=\frac{\xi}{P}.$$

The proof of \eqref{CS} presented here consists analogously with the proof of Theorem $5.2$ in \cite{KirkKlassen90} of two parts. In the first part, we find a path of $\SL(2,\mathbb{C})$ connections on $X$ connecting $\rho$ to an abelian representation $\rho_0$. In fact $\rho_0$ will be an $\SU(2)$ connection on $X.$  In the second part, we then find a path from $\rho_0$ to the trivial representation $\rho^{\text{triv}}$ and we then apply Kirk and Klassens formula \eqref{KK}. The only difference from the proof in \cite{KirkKlassen90} is that we need to explicitly ensure that our paths stay away from parabolic representations. The relevant paths are chosen such that $\lambda,\mu$ are mapped to the maximal $\mathbb{C}^*$ torus of diagonal matrices.

After conjugating by $S_n^{-1}$ we have $\rho(x_n)=\exp\left(\frac{2\pi i l_n}{p_n}\right).$ Consider the subset $$S\subset \text{Hom}\left(\pi_1(Y),\SL(2,\mathbb{C})\right) $$ of representations $\tilde{\rho}$ satisfying
\begin{equation}
\tilde{\rho}(h)=\rho(h), \ [\tilde{\rho}(x_1)]=\left[\exp\left(\frac{\pi il_1}{p_1}\right)\right]\, $$
and
$$\ [\tilde{\rho}(x_j)]=\left[\exp\left(\frac{2\pi il_j}{p_j}\right)\right], \ \text{for }2\leq j\leq n-1.
\end{equation}
where $[Q]$ denotes the $\SL(2,\mathbb{C})$ conjugacy class of $Q \in \SL(2,\mathbb{C}).$ By considering the presentation \eqref{presentation2}, we see that $S$ is naturally homeomorphic to the product of the $n-1$ conjugacy classes
\begin{equation}
S\simeq \left[\exp\left(\frac{\pi il_1}{p_1}\right)\right]\times  \bigtimes_{j=2}^{n-1} \left[\exp\left(\frac{2\pi il_j}{p_j}\right)\right].
\end{equation}
 Therefore the connectedness of $\SL(2,\mathbb{C})$ implies that $S$ is connected. Following a similar argument use in the proof of the previous proposition, we let 
 $$m : S \rightarrow \SL(2,\mathbb{C})$$
  be the algebraic product map. Let $P\subset \SL(2,\mathbb{C})$ be the set of non-diagonalizable elements in $\SL(2,\mathbb{C})$ and observe that $P$ has complex co-dimension one, thus so does $m^{-1}(P) \subset S$, but then it follows that $S' = S- m^{-1}(P)$ is also connected.
 
 Write $\rho=\rho_1$ and observe that $\rho_1 \in S'$. Choose a smooth path $\rho_t$ in $S'$ connecting $\rho_1$ to $\rho_0 \in S'$ given by
\begin{equation}
\rho_0(x_1)=\exp\left(-\frac{\pi i l_1}{p_1}\right), \ \rho_0(x_j)=\exp\left(-\frac{2\pi i l_j}{p_j}\right), \ j=2,...,n-1
\end{equation} and $\rho_0(x_n)=\exp\left(\frac{\pi i l_1}{p_1}+\sum_{j=2}^{n-1}\frac{2\pi i l_j}{p_j}\right).$ By an over all conjugation we can choose the arc $\rho_t$ such that $\rho_t(x_n)=\exp(2\pi if(t))$ for a smooth function $f(t).$  We have $f(0)=\frac{l_1}{2p_1}+\sum_{j=2}^{n-1} \frac{l_j}{p_j}$ and $ f(1)=\frac{l_n}{p_n}.$ Notice that
$f(0)=\frac{\eta}{2}-f(1).$ As $q_n$ is even we have the following two equalities
\begin{align}
\rho_t(\mu)=\rho_t(x_n)^{p_n}\rho_t(h)^{q_n}=&\exp(2\pi i p_nf(t)),
\\ \rho_t(\lambda)=\rho_t(x_n)^{-p_1\cdots p_{n-1}}\rho_t(h)^c &= \exp\left(-2\pi i p_1\cdots p_{n-1}f(t)+vc\pi i\right).
\end{align}
Write $y=vc\in \Z.$ Define $\alpha_1(t)=p_nf(t)$ and $ \beta_1(t)=-\frac{P}{p_n}f(t)+\frac{y}{2}. $ We have that
\begin{align} \label{firstintegral}\begin{split}
-2 \int_0^1 \alpha'_1(t)\beta_1(t) \ \dv t&=-2 \int_0^1 p_nf'(t)\left(-\frac{P}{p_n}f(t)+\frac{y}{2}\right) \ \dv t
\\ &=-2 \int_{f(0}^{f(1)} \left(-Pu+\frac{p_n y}{2}\right) \ \dv u
\\ &= -2\left[-\frac{Pu^2}{2}+\frac{p_ny u}{2}\right]_{u=f(0)}^{u=f(1)}
\\ &=Pf(1)^2-yp_nf(1)-Pf(0)^2+yp_nf(0)
\\ &= Pf(1)^2-Pf(0)^2+yp_nf(0) \text{ mod } \mathbb{Z}.
\end{split}
\end{align}
For the last identity we used that $yp_nf(1)=yl_n \in \Z.$

For the second part we use the fact that $H_1(Y) \simeq \mathbb{Z}$ with generator $\mu$ to conclude that the abelian $\SU(2)$ connection $\rho_0$ can be connected to the trivial representation $\rho^{\text{triv}}$ by a path of $\SU(2)$ representations $\sigma_t$ with $ \sigma_t(\mu)=\exp(2\pi it \alpha_1(0))$ and $ \sigma_t(\lambda)=\exp(2\pi i \beta_1(0)).$ Let $\alpha_0(t)=t\alpha_1(0)$ and $\beta_0(t)=\beta_1(0).$ As $\s(\rho^{\text{triv}})=0,$ we can apply Kirk and Klassen's formula  \eqref{KK} to obtain
\begin{equation}
-\s(\rho) = \s(\rho^{\text{triv}})-\s(\rho)= -2 \int_0^1 \alpha'_0(t)\beta_0(t) \ \dv t-2\int_0^1 \alpha'_1(t)\beta_1(t) \ \dv t.
\end{equation}
We have
\begin{align}
-2\int_0^1 \alpha'_0(t)\beta_0(t) \ \dv t&=-2f(0)(-Pf(0)+\frac{yp_n}{2})=2Pf(0)^2-yp_nf(0).
\end{align}
Comparing this with \eqref{firstintegral} we get that
\begin{align}
-\s(\rho) & = P(f(1)^2+f(0)^2) \ \modu \ \mathbb{Z}
\\ &=P((f(1)+f(0))^2-2f(0)f(1)) \ \modu \ \mathbb{Z}
\\ &= \frac{\xi^2}{4P} \ \modu \ \mathbb{Z}.
\end{align}
For the last equality we used that $2Pf(0)f(1)\in \Z$ and that $$f(0)+f(1)=\frac{\eta}{2}=\frac{\xi}{2P}.$$
This is what we wanted.  \end{proof}

For $x \in \Q$ let $[x]=x \mod  \Z.$ Introduce the set \begin{multline} \mathcal{W}(p_1,...,p_n)= \left\{ \left[ \frac{-m^2}{4P}\right]  : \text{  $m \in \Z$ is divisible by at most $n-3$ of the $p_j$'s}\right\}.
\end{multline}
Recall that the classical complex Chern-Simons values $\CS^*_{\C}(X)$ is the range of the restriction of $\s$ to $\mathcal{M}^*(X,\SL(2,\C))$. Thus we can compute $\CS_{\C}^*(X)$ as a corollary of Proposition \ref{CSvalues}.

\begin{corollary} \label{cor:HUHU} We have that
\begin{equation}
\CS^*_{\C}(X)= \mathcal{W}(p_1,...,p_n).
\end{equation}
\end{corollary}

\begin{proof} It is clear that $\CS_{\C}^*(X) \subset \mathcal{W}(p_1,...,p_n).$ We must show that for any $y \in \mathbb{Z}$ which is not divisible by more than  three of the $p_j$ we can find $l=(l_1,...,l_n) \in L(p_1,...,p_n)$  which solves the congruence equation
\begin{equation} \label{congruence}
y^2 = \left(P\left( \frac{l_1}{p_1}+\sum_{j=2}^{n} \frac{2 l_j}{p_j}\right) \right)^2 \  \modu \ 4P\mathbb{Z}
\end{equation} 
For $x \in \mathbb{Z}$ and $d \in \mathbb{N}$ let $[x]_d$ denote the congruence class of $x$ in the quotient ring $\mathbb{Z} / d\mathbb{Z}.$ Since $p_j$ is odd for $j\geq 2,$ it follows that $4p_1,p_2,...,p_n$ are also pairwise co-prime. Hence the Chinese remainder theorem applies and the natural ring homomorphism $q:\mathbb{Z} \rightarrow \mathbb{Z}/4p_1\mathbb{Z} \oplus_{j=2}^n \mathbb{Z}/p_j\mathbb{Z},$ given by $x \mapsto ([x]_{4p_1},...,[x]_{p_n}),$ descends to an isomorphism of rings 
$$\overline{q}: \mathbb{Z}/4P\mathbb{Z} \overset{\sim}{\rightarrow} \mathbb{Z}/4p_1\mathbb{Z} \ \oplus \ \bigoplus_{j=2}^n \mathbb{Z}/p_j\mathbb{Z} .$$
 It follows that \eqref{congruence} is in fact equivalent to the following $n$ congruence equations
\begin{align}\label{congruence2}
[{y]^2}_{4p_1}&=\left[l_1\prod_{j=2}^n p_j+2\left(\sum_{j=2}^{n} l_j \prod_{t\not=j} p_t\right)\right]^2_{4p_1}, 
\\ [y]^2_{p_j}&=\left[2l_j \prod_{t\not=j}p_t \right]^2_{p_j}, \ \forall j\geq 2.
\end{align} The coprimality conditions ensures that $2\prod_{t\not=j} p_t$ is an invertible element in $\mathbb{Z}/p_j\mathbb{Z}$ and therefore solving the last $n-1$ of the equations in \eqref{congruence2} can indeed be done with $0 \leq l_j\leq (p_j-1)/2.$ 
It remains only to consider the first of the equations in \eqref{congruence2}. To this end we first observe that
$$ \left( l_1 \prod_{j=2}^n p_j + 2 \sum_{j=2}^n l_j  \prod_{t\neq j} p_j\right)^2 = \left( l_1 \prod_{j=2}^n p_j \right)^2 \text{ mod } 4 p_1$$
But then we can solve
$$ y^2 = l_1^2 \prod_{j=2}^n p^2_j \text{ mod } 4 p_1$$
for $0 \leq l_1 \leq 2p_1$. But we also have that
$$ y^2 = (-l_1 \pm 2 p_1)^2 \prod_{j=2}^n p_j^2\text{ mod } 4 p_1.$$
Thus it follows that we can in fact solve  \eqref{congruence2} with $l_1\in \{0,...,p_1\}.$ 
 Thus we have shown that $\CS_{\C}^*(X) = \mathcal{W}(p_1,...,p_n).$  \end{proof}

%
%

Our analysis of the components of the $\SL(2,\C)$ moduli space and the Chern-Simons values now allow us to prove the following.

\begin{theorem} \label{ThmBrieskornSphere} The Chern-Simons action 
$$\s: \pi_0(\mathcal{M}(X,\SL(2,\mathbb{C}))) \rightarrow \mathbb{R}/\mathbb{Z}$$ 
is injective and induces an isomorphism
$$ \pi_0(\mathcal{M}(X,\SL(2,\mathbb{C}))) \cong \mathcal{W}(p_1,\ldots, p_n) \sqcup \{0\}.$$
\end{theorem} 
\begin{proof}
We use the inverse of the isomorphism
$$\overline{q}: \mathbb{Z}/4P\mathbb{Z} \overset{\sim}{\rightarrow} \mathbb{Z}/4p_1\mathbb{Z} \bigoplus_{j=2}^n \mathbb{Z}/p_j\mathbb{Z}$$
to conclude that for each Chern-Simons value in $\mathcal{W}(p_1,\ldots, p_n)$, there is a unique $l\in L(p_1,\ldots, p_n)$ with the given Chern-Simons value, concluding the proof by the last statement of Proposition \ref{rho_l}.

\end{proof}

Thereom \ref{thm:cs}  is a summary of the main results obtained in this section.

\begin{proof}[Proof of Theorem \ref{thm:cs}]
This follows from Theorem 	\ref{thm:decompSL2C}, Corollary \ref{cor:HUHU} and Theorem  \ref{ThmBrieskornSphere}.
\end{proof}

\section{The Borel transform and Complex Cherns-Simons} \label{sec:BOREL}

We now provide the proof of Theorem \ref{2.1}. The reader not familiar with the Borel transform $\mathcal{B}$ and its relation to the Laplace transform is encouraged to read Section \ref{sec:Borel}, before reading the proof of Theorem \ref{2.1}. For a measurable function $g:\R_+ \rightarrow \C$  of sufficient decay, we use the notation $\mathcal{L}_{\R_+}(g)$ for the Laplace transform -- see  Equation  \eqref{eq:Laplacetransform}.

\begin{proof}[Proof of Theorem \ref{2.1}]
We start by giving a characterization of which of the phases in \eqref{phasedecomp} give a non-zero contribution. Introduce for $\mu \in \mathbb{Q}/ \Z$ the set 
\begin{align}\label{R}
\mathcal{T}(\mu) = & \{m=1,...,2P -1 :-m^2/4P=\mu \ \modu \ \mathbb{Z} \}
\\ =& \{m=1,...,2P -1 :g(2\pi i m)=2\pi i \mu \ \modu \ 2\pi i\mathbb{Z} \}
\end{align}

The set of phases $ 2\pi i R^*(X)$ in \eqref{expansionX} consists of the values $g(2\pi i m)=  \frac{-m^2 2\pi i}{4P}$ for which 
\begin{equation} \label{sumresiduesnotzero}
\sum_{x \in \mathcal{T}(-m^2/4P)} \text{Res} \left( \frac{F(y)e^{kg(y)}}{1-e^{-ky}}, y=2\pi i x \right) \not=0,
\end{equation}
for $m=1,...,2P-1$.
Thus, by Corollary \ref{cor:HUHU}, we must prove that if \eqref{sumresiduesnotzero} holds, then there exists $\tilde m \in \mathcal{T}(-m^2/4P)$ such that at most $n-3$ of the $p_j$'s which divide $\tilde m$.

We start by noting that the set of poles of $F$ is given by
\begin{equation} \label{PolesofF}
\mathcal{P}_F=\{ 2\pi im \mid m \in \mathbb{Z} \ \text{and $m$ is divisible by at most $n-3$ of the $p_j$'s} \}.
\end{equation} 
It follows that if $\tilde m$ is divisible by at least $n-2$ of the $p_j,$ then $F(y)$ does not have a pole at $y=2\pi i \tilde m$ and we get for integral $k$
\begin{align}
\text{Res} \left( \frac{F(y)e^{kg(y)}}{1-e^{-ky}}, y=2\pi i \tilde m \right)&=F(2 \pi i \tilde m) e^{kg(2\pi i\tilde m)} \text{Res} \left(\frac{1}{1-e^{-ky}}, y=2\pi i \tilde m \right)
\\ &=F(2 \pi i \tilde m) e^{kg(2\pi i\tilde m)}\frac{1}{k}. \end{align}
As we already noted above, Lawrence and Rozansky checked that all the $k^{-1}$ terms cancels, so it follows that 
\begin{equation}
\sum_{\substack{\tilde m \in \mathcal{T}(-m^2/4P), \\ \text{ and $\tilde{m}$ is divisible by at least $n-2$ of the $p_j$}}}  \text{Res} \left( \frac{F(y)e^{kg(y)}}{1-e^{-ky}}, y=2\pi i \tilde m \right)=0.
\end{equation}
Therefore we see that if \eqref{sumresiduesnotzero} holds, then there is some $\tilde{m} \in \mathcal{T}(-m^2/4P)$ which is divisible by at most $n-3$ of the $p_j.$ 
This establishes  $R^*(X) \subset \CS_{\C}^*$ and we get \eqref{EXPANSIOON}. Observe that as a corollary we obatin for each $\theta \in \CS_{\C}^*$ the formula
\begin{equation} \label{eq:formulaforpolynomial}
e^{2\pi ik \theta} \z_{\theta}(k)=- \sum_{\tilde m \in \mathcal{T}(\theta)}  \text{Res} \left( \frac{F(y)e^{kg(y)}}{1-e^{-ky}}, y=2\pi i \tilde m \right).
\end{equation}

We now turn to $\mathcal{B}(\z_{0}).$ The formal series $\z_{0}$ is the   the asymptotic expansion of the Laplace integral $$\z^{\I}(k)=\frac{1}{2\pi i} \int_{\gamma} F(y) e^{kg(y)} \ \dv y.$$ Let $G$ be the rational function introduced in \eqref{eq:G} and introduce the multivalued function $\mathcal{B}_0(\zeta)$ given by 
\begin{equation}
\mathcal{B}_0(\zeta)= \frac{\kappa}{\pi i \sqrt{\zeta}} G\left( \exp\left(\frac{\kappa\sqrt{\zeta}}{P} \right)\right)= \frac{\kappa i}{4\pi } \frac{\prod_{j=1}^n \sinh\left(\frac{\kappa\sqrt{\zeta}}{p_j}\right)}{ \sqrt{\zeta} \left(\sinh\left(\kappa\sqrt{\zeta}\right)\right)^{n-2}}.
\end{equation} With this notation, the equation for the Borel transform  \eqref{ComputingBorel} which we want to prove, reads as follows
\begin{equation} \label{eq:reduction}
\mathcal{B}(\z_0)= \mathcal{B}_0. 
\end{equation} The function $\mathcal{B}_0$ is related to $F$ as follows
\begin{equation} \label{eq:BORELLO}
\mathcal{B}_0(\zeta)=\frac{\sqrt{ 2 P }}{ \sqrt{\pi i\zeta}} F\left(\sqrt{8\pi i P \zeta}\right).
\end{equation} 
 Now as, $F(-y)=F(y)$ we have a convergent power series expansion valid for $y$ close to $0$
$$F(y)=\sum_{m=1}^{\infty} \frac{F^{(2m)}(0)}{(2m)!} y^{2m}.$$
Therefore Equation \eqref{eq:BORELLO} implies that if we set $$B_m=  \frac{1}{2 \pi i} \sqrt{8\pi i P} \frac{F^{(2m)}(0)}{(2m)!}  \left(8\pi i P \right)^{m}$$ then we have a convergent expansion valid for $\zeta$ close to $0$ of the form
\begin{equation} \label{eq:zetanearzero}
\mathcal{B}_0(\zeta) = \sum_{m=1}^{\infty}  B_m \zeta^{m-1/2}.
\end{equation}
  Introduce the variable $t$ defined by 
  $$-t=g(y)= \frac{ i y^2}{8 \pi P}.$$ 
  Thus 
  $$ \dv y= \frac12 \sqrt{\frac{8 \pi i P }{  t}} \dv t.$$ 
  We now rewrite $\z^{\I}(k)$ as the Laplace transform of $\mathcal{B}_0$ 
\begin{align} \label{ComputeBorel} \begin{split}
\z^{\I}(k)&= \frac{1}{2 \pi i} \int_\gamma F(y) e^{kg(y)} \ \dv y =  \frac{1}{2 \pi i} \int_0^{\infty} e^{-kt} \int_{g=t} \frac{F}{\dv g} \ \dv t
\\ &= \frac{1}{2 \pi i} \int_0^{\infty} e^{-kt} \sqrt{\frac{8 \pi i P}{ t}} F\left(\sqrt{8\pi i Pt}\right) \ \dv t
\\ &= \int_0^{\infty} e^{-kt} \mathcal{B}_0(t) \ \dv t = \mathcal{L}_{\R_+} (\mathcal{B}_0)(k).
 \end{split}
\end{align}
  The existence of the asymptotic expansion \eqref{expansionX}
  
  $$\z^{\I}(k)= \mathcal{L}_{\R_+} (\mathcal{B}_0)(k) \sim_{k \rightarrow \infty}  \z_0(k) $$ can now be obtained by appealing to the first part of Lemma \ref{lem:masterlemma} where we set $B=\mathcal{B}_0$. Here we use the existence of the expansion \eqref{eq:zetanearzero}. Therefore the desired identity  \eqref{ComputingBorel}  $\mathcal{B}(\z_0)=\mathcal{B}_0$ follows from the second part of Lemma \ref{lem:masterlemma} and the convergence of the expansion \eqref{eq:zetanearzero}. 

  As $F(-y)=F(y)$ we note that the factor
  
   $$\zeta \mapsto F\left(\sqrt{8\pi iP\zeta}\right)$$
    gives a well-defined meromorphic function. Thus $\mathcal{B}(\z_{0})(\zeta)$ is a multivalued meromorphic function with a square root singularity at $0$ and with  singularities for $\sqrt{8\pi iP\zeta} \in  \mathcal{P}_F$ where $\mathcal{P}_F$ is the set of poles of $F(y).$ This set was computed above (see equation \eqref{PolesofF}) and we conclude that the poles of  $\mathcal{B}(\z_{0})(\zeta)$ occur at
    
$$\zeta_m=\frac{ -\pi  m^2}{2iP}=\frac{- m^2}{4P} \frac{2\pi }{i}$$ with $m\in \mathbb{Z}$ being divisible by less than or equal to $n-3$ of the $p_j$'s.  This concludes the proof of \eqref{Inclu}.
 \end{proof}
 
  It is of course expected that only a Chern-Simons invariant $\theta$ of a flat $\SU(2)$ connection have a non-vanishing polynomial $\z_{\theta} \not=0,$ i.e. $$R^*(X)=\s(\mathcal{M}^*\left(X,\SU(2)\right) ).$$

\subsection{Resummation of the WRT invariant}
We now turn to the resummation of the normalized WRT invariant $\tilde{\z}_k(X)$. Recall that for $\mu \in \mathbb{Q}/ \Z$ we introduced the set
\begin{align}\label{R}
\mathcal{T}(\mu) = & \{m=1,...,2P -1 :-m^2/4P=\mu \ \modu \ \mathbb{Z} \}.
\end{align}
We also introdude the residue operator $\mathcal{I}_\mu$ which for a meromorphic function $\hat{\varphi}$ is given by
\begin{align}  \label{Lmu}  \mathcal{I}_{\mu}(\hat{\varphi})(\xi) = -& \sum_{x \in \mathcal{T}(\mu)} \Res \left (\frac{\exp\left(\xi \frac{iy^2}{8\pi P}\right)}{(1-e^{-\xi y})} \frac{y}{ 4P}\hat{\varphi}\left(\frac{y^2}{i8\pi P }\right) , y=2\pi i x\right)\end{align}
Observe that by definition $\mathcal{T}(\mu)$ is empty for all but finitely many $\mu \in \mathbb{Q} / \mathbb{Z}$ and therefore $\mathcal{I}_{\mu}$ is $0$ for all but these finitely many $\mu.$

\begin{corollary}  \label{2.2}
	The polynomials $\z_{\theta}$ and the quantum invariant $\widetilde{\z}_k(X)$  are determined by $\mathcal{B}(\z_{0})$  as follows
	\begin{align} \label{Poooly}  \z_{\theta} (k) & = e^{-2\pi i k \theta} \mathcal{I}_{\theta}\left(\mathcal{B}(\z_{0})\right)(k).
	\\   \widetilde{\z}_k(X) &= \mathcal{L}_{\R_+}( \mathcal{B}(\z_{0}))(k)+\sum_{\theta \in \frac{i}{2\pi} \Omega \ \modu \ \Z} \mathcal{I}_{\theta}\left(\mathcal{B}(\z_{0})\right)(k) \label{Resuuumation} .
	\end{align}
\end{corollary} 
The identity \eqref{Resuuumation} of Corollary \ref{2.2} is reminiscent of the typical resummation process from resurgence \cite{Schiappa18,Dorigoni14}. The Ohtsuki-series is known to determine $\tau_k(X).$ The new insight provided by resurgence is that it does so via resummation as stated in Corollary \ref{2.2}.

We now prove Corollary \ref{2.2}. 
\begin{proof}
	It easily follows from \eqref{ComputingBorel} that
	\begin{equation} \label{TYTY}
	F(\zeta)= \mathcal{B}(\z_{0})\left(\frac{\zeta^2}{i8P\pi}\right) \frac{\zeta }{4P}.
	\end{equation} Recall from  \eqref{eq:formulaforpolynomial} that \begin{equation} 
	e^{2\pi ik \theta} \z_{\theta}(k)=- \sum_{\tilde m \in \mathcal{T}(\theta)}  \text{Res} \left( \frac{F(y)e^{kg(y)}}{1-e^{-ky}}, y=2\pi i \tilde m \right).
	\end{equation}
From this and Equation \eqref{TYTY} one easily obtains \eqref{Poooly}. 

In the proof of Theorem \ref{2.1} we obtained the following exact formula
\begin{equation} \label{eq:DecompWRT}
\tilde{\z}_k(X) = \mathcal{L}_{\R_+}(\mathcal{B})(k)+ \sum_{\theta \in \CS_{\C}^*(X)} e^{2\pi i k \theta} \z_\theta(k).
\end{equation}
Thus we see that \eqref{Resuuumation} follows from this formula and \eqref{Poooly}. 
\end{proof}

\subsection{Resurgence of the generating function} \label{sec:GeneratingFunction}

Let $G$ be a simple, simply connected compact Lie group, and let $\tau_{G,k}$ be the level $k$ Reshetikhin-Turaev TQFT constructed from the quantum group $U_q(\mathfrak{g})$, where $\mathfrak{g}$ is the complexification of the Lie algebra $\mathfrak{h}$ of $G$. Let $\check{h}$ be the dual Coxeter number of $\mathfrak{h}$, and set $\check{k}=k+h$. For a closed oriented three manifold $M$ (possibly containing a colored framed link) we consider the normalized invariant \begin{equation}
\z_{G,k}(M)= \frac{\tau_{G,\check{k}}(M)}{\tau_{G,\check{k} }(\Sph^2 \times \Sph^1)}.
\end{equation}
Let $z$ be a formal variable and consider the generating function $$\z_{G}(M;z) \in \C[[z]]$$
given by
\begin{equation} \label{eq:GeneratingFunction}
\z_G(M;z) = \sum_{k=0}^{\infty} \z_{G,k}(M) z^k.
\end{equation}
By work of Garoufalidis  $\z_{G}(M;z) $ is known to be convergent on the unit disc. Motivated by the paradigm of analytic continuation and resurgence, Garoufalidis posed the following conjecture

\begin{conjecture}[\cite{Garoufalidis08}] \label{conj:GeneratingFunction}
	The generating function $\z_{G}(M;z)$ has an analytic continuation to $\C\setminus e \Lambda$ where $e \Lambda$ is a finite set containing zero and the exponentials of the negatives of the complex classical Chern-Simons values.
\end{conjecture}

In other words, the conjecture is that the generating function $\z_{G}(M;z)$ determines the germ at zero of a resurgent function. This conjecture is formally motivated from resurgence of Laplace integrals and the (non-rigorous) path integral formula for the WRT invariant, as explained in \cite{Garoufalidis08}.

We now specialize to the case of the Seifert fibered homology sphere $X$ and $G=\SU(2)$. Set $K=k+2$ and consider the generating function for the normalized quantum invariant $\tilde{\z}_k(X)$ given by
\begin{equation} \label{eq:NGeneratingFunction}
\tilde{\z}(X;z) =\sum_{k=0}^{\infty} \tilde{\z}_K(X) z^k  \in \C[[z]].
\end{equation}

For $s \in \C$ consider the polylogarithm
\begin{equation} \label{eq:polylogarithm}
\Li_{s}(z) = \sum_{n=1}^{\infty} \frac{z^n}{n^s}.
\end{equation}
For $s=-m, m \in \N$ the polylogarithm is exact and in fact a rational function 
\begin{equation}
\Li_{-m} (z)= \left( z\frac{\partial }{\partial z} \right)^m \left( \frac{z}{1-z}\right).
\end{equation}
We introduce the following notation for the exponentials of the negatives of the classical complex Chern-Simons values
\begin{equation}
e \Lambda = \exp\left(- 2\pi i \CS^*_{\C}(X) \right).
\end{equation}
We prove the following proposition.

\begin{proposition} \label{thm:ExactGeneratingFunction} The generating function $\tilde{\z}(X;z)$ is the germ at zero of a holomorphic function
	$\tilde{\z} \in \mathcal{O}(\C \setminus e\Lambda )$ given by the following formula
	
	\begin{align} \begin{split} \label{eq:exact}
	\tilde{\z}(X;z)&= \int_0^{\infty} \frac{e^{-2y} \mathcal{B}(\z_0)(y) }{1-ze^{-y}} \dv y
	\\&+\sum_{\theta \in \CS_{\C}^*(X)}  e^{4\pi i \theta} \sum_{j=0}^{n-3} \frac{\z_{\theta}^{(j)}(0)}{j!} \left(2^j +\sum_{l=0}^j 2^{j-l} \binom{j}{l} \Li_{-l} \left(ze^{2\pi i \theta}\right)\right). 
	\end{split} 
	\end{align}
\end{proposition}

\begin{proof}
	From Equation \eqref{eq:DecompWRT} it follows that \begin{align}  \label{eq:kl} \begin{split} 
	\tilde{\z}(X;z)=& \sum_{k=0}^{\infty} z^k \int_0^{\infty} e^{-y(k+2)} \mathcal{B}(\z_0)(y) \dv y 
	\\ & + \sum_{ \theta \in \CS_{\C}(X)}  \sum_{k=0}^{\infty}  z^k \z_{\theta}(k+2) e^{2\pi i (k+2) \theta} . \end{split}
	\end{align}
	The first term can be simplified by interchanging summation and integration and then using the geometric series expansion
	
	\begin{align} \label{eq:jiji}  \begin{split} 
	\sum_{k=0}^{\infty} z^k \int_0^{\infty} e^{-y(k+2)} \mathcal{B}(\z_0)(y) \dv y &= \int_0^{\infty} \sum_{k=0}^{\infty}  (ze^{-y})^k \mathcal{B}(\z_0)(y) e^{-2y} \ \dv y
	\\ &= \int_0^{\infty} \frac{e^{-2y}\mathcal{B}(\z_0)(y) }{1-ze^{-y}} \ \dv y.   \end{split} 
	\end{align}
	This can be justified by standard complex analysis arguments. To complete the proof, we can consider seperately each term in \eqref{eq:kl} corresponding to a complex Chern-Simons value $\theta \in \CS_{\C}(X)$. We get that
	\begin{align} \label{eq:uiui} \begin{split} 
	&	\sum_{k=0}^{\infty} z^k \z_{\theta}(k+2) e^{2\pi i (k+2)\theta}
	\\ &=e^{4\pi i \theta} \sum_{j=0}^{n-3} \frac{\z_{\theta}^{(j)}(0)}{j!}\sum_{k=0}^{\infty} (k+2)^j (e^{2\pi i \theta }z)^k 
	\\ &= e^{4\pi i \theta} \sum_{j=0}^{n-3} \frac{\z_{\theta}^{(j)}(0)}{j!}\left(2^j +\sum_{l=0}^{j} 2^{j-l}\binom{j}{l} \sum_{k=1}^{\infty} k^l (e^{2\pi i\theta } z)^k \right)
	\\ &=e^{4\pi i \theta} \sum_{j=0}^{n-3} \frac{\z_{\theta}^{(j)}(0)}{j!} \left(2^j +\sum_{l=0}^{j} 2^{j-l}\binom{j}{l} \Li_{-l}(e^{2\pi i\theta}z) \right). \end{split}
	\end{align}
	In the last equality, we used the series expansion \eqref{eq:polylogarithm} of the polylogarithm. By substituting the identities \eqref{eq:jiji} and \eqref{eq:uiui} into \eqref{eq:kl} we obtain the desired identity \eqref{eq:exact}.
\end{proof}

  \section{A resurgence formula for the GPPV invariant }  \label{sec:Zhat}
  
  We now turn to the $q$-series invariant $\hat{\z}_0(X,q).$  We follow \cite{GukovCiprian19}. Let $(\Gamma,m)$ be an ordered weighted tree, i.e. $\Gamma$ is a tree together with an ordering of its set of vertices $V$ and $m$ is a map $m:V\rightarrow \Z.$  Set $s=\lvert V \rvert$ and let $M=M(\Gamma,m)$ be the $s\times s$ matrix with entries given by
  $$M_{i,j}= \begin{cases} m_v & \text{if} \ v_i=v_j=v,
  \\ 1 & \text{if}  \ v_i \text{ and }  v_j \ \text{are joined by an edge,}
  \\ 0 & \text{otherwise.} \end{cases}$$
  We say $M$ is weakly negative definite if $M$ is invertible and $M^{-1}$ is negative definite on the subspace of $\Z^s$ spanned by vertices of degree at most $3.$ Let $Y=Y(\Gamma,m)$ be the oriented three manifold with surgery data $L=L(\Gamma,m)$ constructed as follows. For each vertex $v$ the link $L$ has an unknotted component $L_v$ with framing $m_v,$ and $L_v$ is chained together with $L_w$ if and only if $v$ and $w$ are joined by an edge. We call $Y$ a plumbed manifold with plumbing graph $\Gamma.$
  
  We recall that two plumbed three manifolds $Y$ and $Y'$ are diffeomorphic if and only their plumbing graphs are related by Neumann moves. 
  
  When $Y$ is a plumbed manifold with weakly negative definite plumbing graph and $Y$ is not necessarily a homology $3$-sphere, the $q$-series invariant $\hat{\z}_a(Y,q)$ depend on a label $a,$ whose precise  meaning is subtle. Originally, these labels where thought to be abelian flat connections, later $\Spin$ structures, and for mapping tori of genus one mapping classes one has include "almost abelian" flat connections (see \cite{2019arXiv191108456C}). As $X=\Sigma(p_1,...,p_n)$ is a homology three sphere, we have $a=0,$ and need not go deeper into this discussion. For the sake of completeness however, we recall the GPPV-formula definition as it is stated in terms of $\Spin$-structures. First we recall how $\Spin$-structures can be described in terms of the adjacency matrix $M.$ This is thorougly explained in \cite{GukovCiprian19}. Let $Y$ be a plumbed three manifold with plumbing graph $\Gamma.$ Let $s=\lvert V \rvert.$ Let $\vec{m} \in \Z^{s}$ be the weight vector, i.e. $m_j=m(v_j).$ Let $\vec{\delta} \in \Z^s,$  be the degree vector i.e. $\delta_j=\text{deg}(v_j).$ We have isomorphisms $$
  \Spin(Y) \simeq (\Z^s+\vec{m})/ 2M\Z^s \simeq (\Z^s+\vec{\delta})/ 2M\Z^s.$$
  These isomorphisms  are compatible with Neumann moves as explained in \cite{GukovCiprian19}. We now recall the GPPV-formula \eqref{eq:GPPV}.
  \begin{definition}[\cite{GukovPeiPutrovVada18}] \label{df:GPPV}
  	Let $Y$ be a plumbed three manifold with weakly definite plumbing graph $\Gamma.$ Let $\phi$ denote the number of positive eigenvalues of $M$ and let $\sigma$ denote the signature of $M.$ Let $a=[ \vec{a}] \in (\Z^s+\vec{\delta})/ 2M\Z^s \simeq \Spin(Y).$ The $\hat{Z}$-invariant of $(Y,a)$ is given by 
  	\begin{equation} \label{eq:GPPV}
  	\hat{\z}_a(Y;q)=(-1)^{\phi}q^{\frac{3\sigma-\sum_v m_v}{4}} \cdot v.p. \oint_{\lvert z_v \rvert=1} \prod_{v\in V} \frac{\dv z_v}{2\pi iz_v}\left(z_v-\frac{1}{z_v}\right)^{2-\text{deg}(v)} \Theta_a^{-M}(\vec{z})
  	\end{equation}
  	where $v.p.$ denotes the principal value and
  	\begin{equation} \Theta_a^{-M}(\vec{z})=\sum_{\vec{l}\in 2M\Z^s+\vec{a}}
  	q^{-\frac{(\vec{l},M^{-1}\vec{l}))}{4}} \prod_{v \in V} z_v^{l_v}
  	\end{equation}
  \end{definition}
  \begin{remark} The invariance of \eqref{eq:GPPV} under Neumann moves is proved by Gukov and Manulescu in \cite{GukovCiprian19}. \end{remark} We recall that the principal value $v.p.$ is defined such that for every sufficiently small $\epsilon >0$ we have $$v.p. \oint_{\lvert z \rvert =1}= \frac{1}{2}\left(\oint_{\lvert z \rvert =1+\epsilon}+\oint_{\lvert z \rvert =1-\epsilon}\right).$$  
  
  \subsection{Proof of Theorem \ref{thm:zedhat}}
\begin{figure}[htp]
  	\centering
  	\includegraphics[scale=0.8]{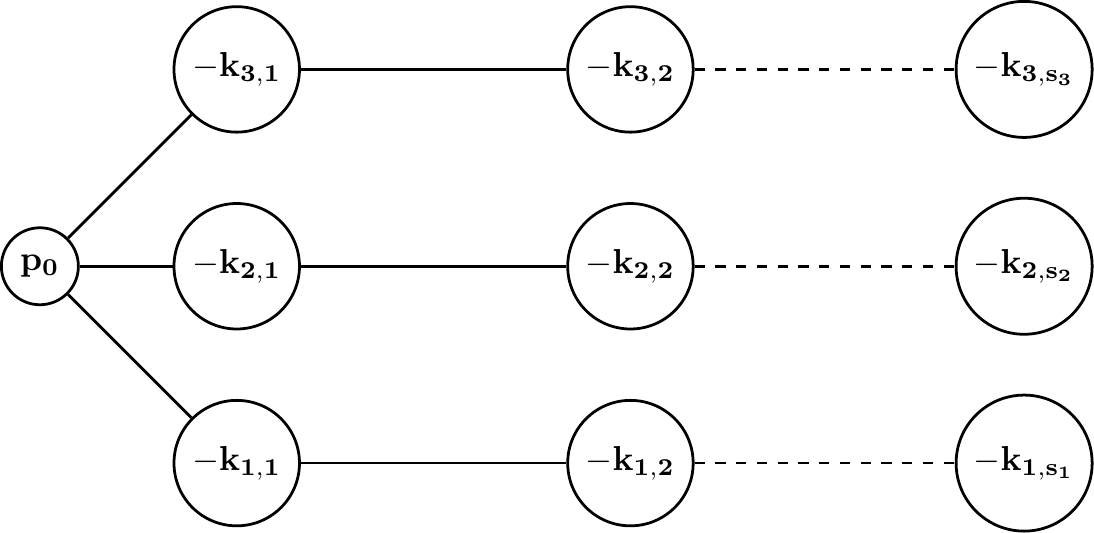}
  	\caption{Plumbing graph for $X$ in the case $n=3.$}
  	\label{fig:graph}
  \end{figure}

  We now consider $X= \Sigma(p_1,...,p_n)$ in more detail. Choose $q_1,..., q_n \in \N$ such that for each $j=1,...,n$ we have $(p_j,q_j)=1, $ $p_j \leq q_j$ and
  \begin{equation}
  p_0:=\frac{-1}{P}-\sum_{j=1} \frac{q_j}{p_j} \in \Z_{<0}.
  \end{equation}
  Then $e=-1/P<0$ is the  Seifert Euler number. Choose a continued fraction expansion of $q_j/p_j$ for each $j=1,...,n$
  $$\frac{q_j}{p_j}=k_{j,1}- \frac{1}{k_{j,2}- \frac{1}{\ddots-k_{j,s_j}}}.$$ As explained in \cite{GukovCiprian19}, $X$ has a negative definite plumbing graph $\Gamma$ defined as follows. The graph $\Gamma$ is star-shaped with $n$ arms and  central vertex $v_0$ with weight $p_0.$ For each $j=1,...,n$ the $j'$th arm has $s_{j}$ vertices. If these are ordered $(v_{j,1},...,v_{j,s_j})$ with $v_{j,1}$ being closest to the central vertex $v_0$, then $v_{i,j}$ have weight $-k_{j,i}.$ This graph is illustrated for $n=3$ in Figure \ref{fig:graph} 

  Before proving Theorem \ref{thm:zedhat}, we first give a formula for the rational exponent $\Delta \in \mathbb{Q}.$  For each $j=1,...,n$ let $X_j$ be the plumbed manifold whose graph $\Gamma_j$ is identical to $\Gamma$ except that on the $j'$th arm, we delete the terminal vertex $v_{j,s_j}.$  Define $h_j \in \N$  as
  $$h_j= \lvert \h_1(X_j,\Z) \rvert. $$
  Observe that the total number of vertices of $\Gamma$ is given by$s=1+\sum_{j=1}^n s_j.$ Define $\Delta \in \mathbb{Q}$ by
  \begin{equation} \label{def:delta}
  \Delta=-\frac{1}{4}\left(\sum_{j=1}^n h_j- 3s- p_0+ \sum_{ j=1}^n \left(-\frac{P}{{p_j}^2}+ \sum_{i=1}^{s_j} k_{j,i} \right) \right)
  \end{equation}

  We now prove Theorem \ref{thm:zedhat}.
  
  \begin{proof}
  	Recall that  $q=\exp(2\pi i \tau)$  where $\tau \in \mathfrak{h}$.  For the sake of notational simplicity, we also introduce the paramter $h=2\pi i \tau$ so that $q=\exp(h)$.  We start by proving that
  	\begin{equation} \label{eq:integralformula}
  	\I(h)=\Psi(q),
  	\end{equation}
  	where $\Psi(q)$ is the series introduced in \eqref{eq:q-series} and $\I(h)$ is the contour integral introduced in \eqref{eq:contourintegral} (with $h=2\pi i \tau)$. Observe that for the purpose of proving \eqref{eq:integralformula} we can and will assume that
  	$$\tau \in i \R_{>0},$$
  	because if the identity \eqref{eq:integralformula} holds true on this half-line, it has to hold on the entire upper half-plane $\mathfrak{h}$, since both functions are holomorphic in $\mathfrak{h}$.

  	Set  \begin{equation} \label{def:Bnorm}
  	\mathcal{B}(\zeta)=\frac{1}{2\sqrt{P}}\mathcal{B}(\z_{0})\left(\frac{\zeta}{2\pi i}\right).
  	\end{equation}
  	For all $t\in \C$ with $\sqrt{t} \in \{z \in \C\mid \re(z)<0\}$ the normalized Borel transform $\mathcal{B}$ satisfies by Theorem \ref{2.1}
  	\begin{equation} \label{eq:twistedBorel} \mathcal{B}(t)= \frac{1}{\sqrt{t}} \sum_{m=m_0}^{\infty} c_m \exp\left(m\sqrt{\frac{t}{P}} \right)= \frac{1}{\sqrt{t}}  G\left(\exp\left(\sqrt{\frac{t}{P}} \right)\right)\end{equation}
  	where for all $m\in \N$ we have that
  	$$c_m=(-1)^n\chi_m.$$
  	For each $m\in \Z_{\geq m_0}$ introduce the polynomial
  	$$p_m(w)=-\frac{w^2}{h}+\frac{m}{\sqrt{P}}w.$$
  	This is a Morse function with a unique saddle point at 
  	$w_m=\frac{h m}{2 \sqrt{P}}$
  	and we have that
  	$$p_m(w_m)=h \frac{m^2}{4P}.$$
  	Let $\Disc(R)$ be the closed ball centered at the origin with radius $R> \lvert m_0 \rvert$. We can deform $i\R$ slightly to a contour $\Delta_m \subset \{z\in \C: \re(z)<0\}\cup \Disc(R),$ which passes through the saddle point $w_m$ and such that the function given by  $$w \mapsto \exp(p_m(w)),$$ have exponential decay along $\Delta_m.$ The orientation of $\Delta_m$ is as depicted in Figure \ref{fig:Delta}. We remark that Figure \ref{fig:Delta} depicts the situation where $m_0\geq 0$. 
  	\begin{figure}[htp]
  		\centering
  		\includegraphics[scale=0.8]{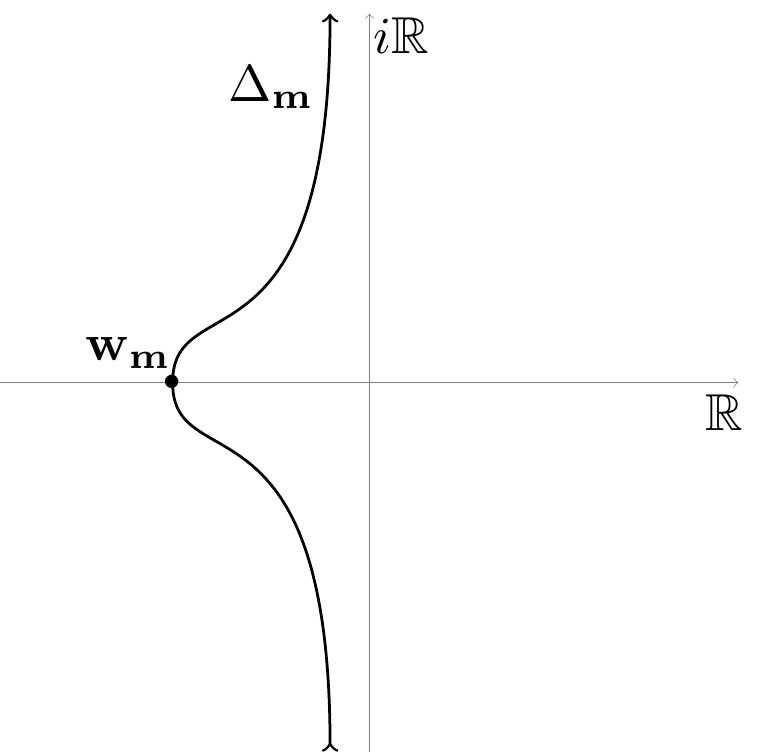}
  		\caption{The contour $\Delta_m.$}
  		\label{fig:Delta}
  	\end{figure}
 Recall that if $\Gamma$ is a steepest descent contour through the unique saddle point of a degree two polynomial $p(z)=-\alpha z^2+\beta z,$ then we have the following exact formula known as Gaussian integration
  	$$\int_{\Gamma} \exp(p(z)) \ \dv z= \sqrt{\frac{\pi}{\alpha}}\exp\left(\frac{\beta^2}{4\alpha}\right).$$  Applying Gaussian integration to the polymonials $p_m$ gives us  the following identity
  	\begin{equation} \label{eq:psipsipsi}
  	(-1)^n\Psi(q)=\sum_{m=m_0}^\infty c_m q^{\frac{m^2}{4P}}=\sum_{m=m_0}^\infty c_m \frac{1}{\sqrt{\pi h}} \int_{i\R}  \exp\left(p_m(w) \right)\dv w.
  	\end{equation}
  	Choose a small postive parameter $\delta>0$ and introduce the contour $$\Delta_0 =e^{i\delta}i\R_+ \cup e^{-i\delta}i\R_{-} \subset \{z \in \C: \re(z)<0\}.$$ 
  	Let $\Upsilon \subset  \{z\in \C: \Re(z)<0\}$ be the Hankel contour which encloses $\R_{-}$ and satisfies $\sqrt{\Upsilon}=\Delta_0.$ The orientation of these contours are given as in Figure \ref{fig:UpsilonDelta}. \begin{figure}[htp]
  		\centering
  		\includegraphics[scale=0.8]{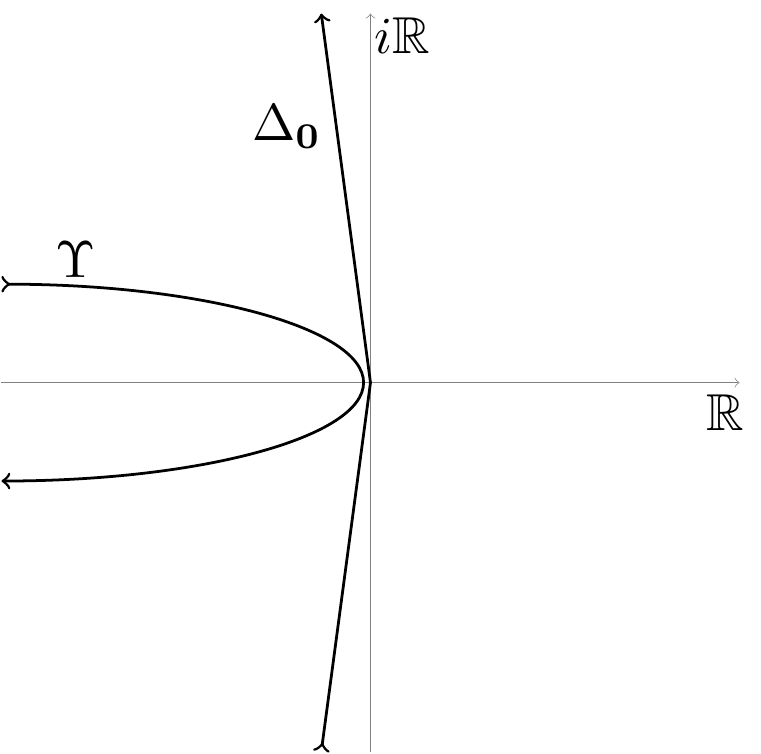}
  		\caption{The contours $\Upsilon$ and $\Delta_0.$}
  		\label{fig:UpsilonDelta}
  	\end{figure}   We let 
  	$$\Upsilon_\pm = \Upsilon \cap \{z \in \C: \im(z) \in \R_{\mp }\},$$ 
  	so that $\sqrt{\Upsilon_+}= e^{i\delta}i\R_+ $ and $-\sqrt{\Upsilon_-}= e^{-i\delta}i\R_{-}$ where $\sqrt{\cdot}$ denotes the principal branch of the square root. Introduce the variable $$w^2=v.$$ As $\Delta_0$ is a small deformation of $\Delta_m$ for each $m\in \Z_{\geq m_0},$ we obtain
  	\begin{align} \begin{split} \label{eq:hammerdammer}
  	\sum_{m=m_0}^\infty  \frac{c_m}{\sqrt{\pi h}} \int_{i\R} \exp(p_m(w)) \ \dv w &= \sum_{m=m_0}^\infty   \sum_{\epsilon =\pm 1 } \frac{c_m}{\sqrt{\pi h}}  \ \int_{\Upsilon(\epsilon)}  \frac{\exp\left(-\frac{v}{h}+\epsilon \frac{m}{\sqrt{P}}\sqrt{v} \right)}{\epsilon 2\sqrt{v}} \ \dv v
  	\\ &= \sum_{\epsilon=\pm 1}\frac{1}{2\sqrt{\pi h}}\int_{\epsilon \Upsilon(\epsilon)} \frac{\exp(-\frac{v}{h})}{\sqrt{v}} \sum_{m=m_0}^{\infty} c_m e^{\epsilon \frac{m\sqrt{v}}{\sqrt{P}}} \ \dv v
  	\\ &= \sum_{\epsilon=\pm 1 } \frac{1}{2\sqrt{\pi h}}\int_{\epsilon \Upsilon(\epsilon)} \exp\left(-\frac{v}{h}\right) \mathcal{B}(v) \ \dv v
  	\end{split}
  	\end{align}
  	In the second equality of \eqref{eq:hammerdammer} we used that $\epsilon \sqrt{v} \in \{z \in \C: \re{z}<0 \}$ for all $v \in \Upsilon(\epsilon),$ and the contour $-\Upsilon(-1)$ denotes $\Upsilon(-1)$ but oriented in the direction from the origin and towards infinity.  In the third equality of \eqref{eq:hammerdammer} we used equation \eqref{eq:twistedBorel} and the identity
  	\begin{equation}
  	G(z)=G\left( \frac{1}{z} \right),
  	\end{equation}
  	which follows directly from the definition of $G$.
  	Now introduce the variable
  	$$\xi=\frac{v}{2\pi i}.$$
  	This identifies (up to a small deformation) the $v$ contour $\Upsilon_+- \Upsilon_-$  with the $\xi$  contour $\varGamma$ introduced in Figure \ref{fig:ContourGamma}. 

  We have 
  	\begin{multline} \label{eq:ey}
  	\sum_{\epsilon=\pm 1 } \frac{1}{2\sqrt{\pi h}}\int_{\epsilon \Upsilon(\epsilon)} \exp\left(-\frac{v}{h}\right) \mathcal{B}(v) \ \dv v
  	\\= \sum_{\epsilon=\pm 1 } \frac{2\pi i }{2\sqrt{\pi h}}\int_{\exp(\epsilon \frac{\pi i}{4})i\R_+}  \exp\left(-\frac{2\pi i \xi}{h}\right) \mathcal{B}(2 \pi i \xi) \ \dv \xi
  	\\ = \sum_{\epsilon=\pm 1 } \frac{2\pi i }{2\sqrt{\pi h}} \frac{\sqrt{2\pi i}}{8\kappa} \int_{\exp(\epsilon \frac{\pi i}{4})i\R_+}  \exp\left(-\frac{2\pi i \xi}{h}\right)\mathcal{B}(\z_{0})(\xi) \ \dv \xi.
  	\end{multline}
  	In the last equality of  \eqref{eq:ey} we used Equation \eqref{def:Bnorm}, which relates $\mathcal{B}$ and $\mathcal{B}(\z_0).$ Now recall that $\kappa=\sqrt{2\pi i P}$, since $H=1$ and recognize the pre-factor in the last line of \eqref{eq:ey} as
  	\begin{equation}
  	\frac{2\pi i }{2\sqrt{\pi h}} \frac{1}{\sqrt{P}2} =(-1)^n\frac{\lambda}{\sqrt{\tau}},
  	\end{equation}
  	where $\lambda$ is the scalar introduced in the statement of Theorem \ref{thm:zedhat}. By combining
  	Equations \eqref{eq:psipsipsi}, \eqref{eq:hammerdammer} and \eqref{eq:ey}, we see that Equation \eqref{eq:integralformula} holds. 
  	
  	Write $\hat{\z}_0(X;q) =\hat{\z}_0(q)$. We now  show that 
  	\begin{equation} \label{eq:numerodos}
  	\Psi(q)=q^{\Delta}\hat{\z}_0(q),
  	\end{equation}
  	where $\Delta \in \mathbb{Q}$ is the scalar introduced in \eqref{def:delta}. This will establish \eqref{eq:zedhat} and thereby finish the proof. We start with $\hat{\z}_0(q).$ By Definition \ref{df:GPPV} and since in this case $\phi =0$, we have that
  	\begin{equation} \label{eq:G1}
  	\hat{\z}_0(q)=q^{\frac{3\sigma-\sum_v m_v}{4}}\sum_{\vec{l}\in 2M\Z^s}
  	q^{-\frac{(\vec{l},M^{-1}\vec(l))}{4}}  \oint_{\lvert z_v \rvert=1} \prod_{v\in V} \frac{\dv z_v}{2\pi iz_v}\left(z_v-\frac{1}{z_v}\right)^{2-\text{deg}(v)} z_v^{l_v}.
  	\end{equation}
  	Here it is understood that we have taken the principal value of the integral as explained above. Recall that for a Laurent series $a(z)=\sum_{j \in \Z} a_jz^j$ we have that
  	$$v.p. \oint_{\lvert z \rvert =1} \frac{\dv z}{2\pi i z} a(z)=a_0.$$
  	For our star-shaped plumbing graph $\Gamma,$ the  non-zero  contributions to  \eqref{eq:G1} comes from 
  	$\vec{l}\in 2M\Z^s$ with
  	$l_w=0$ for all of the entries corresponding to an  internal vertex $w$ of an arm, and $l_v=\pm 1$ if $v$ is a terminal vertex of an arm and then from the central vertex $v_0$, which we will now consider. 
  	
  	In comparing $\hat{\z}_0(q)$ with $\Psi(q)$ it is useful to introduce the integer sequence $\{a_j\}_{j=0}^{\infty}$ determined by
  	\begin{equation} \label{eq:a_j}
  	\left(t-t^{-1}\right)^{2-n}= \begin{cases} 
  	\sum_{j=0}^{\infty} a_j t^{2j+n-2}, \ & \text{if} \ \lvert t \rvert <1
  	\\ 
  	\\
  	\sum_{j=0}^{\infty} (-1)^n a_j t^{-2j-n+2}, \ & \text{if} \ \lvert t \rvert >1.
  	\end{cases}
  	\end{equation}  
  The $a_j$'s can be explicitly evaulated: By the formula for the geometric series $(1-t)^{-1}= \sum_{j=0}^{\infty} t^j$ and Cauchy multiplication of power  series, we see that
  $$\frac{1}{(1-t)^m}= \sum_{j=0}^{\infty} \left( \sum_{ j_1+\cdots+ j_m=j} 1\right) t^j=  \sum_{j=0}^{\infty} {{j+m-1}\choose{j}}  t^j,$$
  and therefore one sees that
  $$a_j=(-1)^n {{j+n-3}\choose{j}}.$$
  However for the comparison of $\hat{\z}_0(q)$ and $\Psi(q)$  given below, we don't need the closed form for $a_j$, but rather the equation \eqref{eq:a_j}.
  	
  	Write $z_{v_0}=z$ and $l=l_{v_0}$. We obtain
  	\begin{equation}
  	v.p. \oint_{\lvert z \rvert =1} \frac{\dv z}{2\pi i z} \left(z-z^{-1}\right)^{2-n} z^{l}= \begin{cases}
  	\frac{1}{2}a_{\frac{l-n+2}{2}} \ & \text{if} \ 2-n-l \in 2 {\mathbb Z}_+,
  	\\
  	\\ \frac{1}{2}(-1)^na_{\frac{2-n-l}{2}} & \ \text{if} \ l-n+2 \in 2{\mathbb Z}_+,
  	\\
  	\\ a_0 \ & \text{if}\  l=0, n= 2.
  	
  	\end{cases}
  	\end{equation}

  	We know that the adjacency matrix $M$ is unimodular, and so $M\Z^s=\Z^s.$ Define a map $\vec{l}:\{\pm 1 \}\times  \N \times \{\pm 1\}^{n} \rightarrow \Z^s$  as follows: For the central vertex $v_0$ we have
  	\begin{equation}\vec{l}(\varepsilon,j,\epsilon)_{v_0}=\varepsilon(-2j+2-n).
  	\end{equation}
  	For $m=1,...,n$ and the terminal vertex $v$ of the $m'$th arm, we have
  	$$\vec{l}(\varepsilon, j,\epsilon)_{v}=\epsilon_m,$$
  	and for every internal vertex $w$ of the arms, we have
  	$$\vec{l}(\varepsilon,j,\epsilon)_{w}=0.$$
  	With this notation, the above considerations show that
  	\begin{equation}
  	\hat{\z}_0(q)=q^{-\frac{3s+\sum_v m_v}{4}} \sum_{\varepsilon=\pm 1} \sum_{r=0}^{\infty} \sum_{\epsilon \in  \{\pm 1\}^{n}} (-1)^{\frac{(1-\varepsilon)(n-2)}{2}} \frac{a_r}{2} \left(\prod_{j=1}^n \epsilon_j \right)
  	q^{-\frac{\left\langle \vec{l}(\varepsilon ,r,\epsilon),M^{-1} \vec{l}(\varepsilon,r,\epsilon), \right\rangle}{4}}.
  	\end{equation}
  	If we apply the symmetry that simultaneously changes the sign of all $\epsilon_j$ and $\varepsilon,$ then we we obtain
  	\begin{equation}
  	\hat{\z}_0(q)=(-1)^{n}q^{-\frac{3s+\sum_v m_v}{4}} \sum_{r=0}^{\infty} \sum_{\epsilon \in  \{\pm 1\}^{n}} a_r  \left(\prod_{j=1}^n \epsilon_j \right) q^{-\frac{\left\langle \vec{l}(-1,r,\epsilon),M^{-1} \vec{l}(-1,r,\epsilon), \right\rangle}{4}}.
  	\end{equation}
  	The quadratic form $$
  	\vec{l} \mapsto \langle \vec{l} ,M^{-1}\vec{l}\rangle/4$$ was computed for $n=3$ in \cite{GukovCiprian19} in their proof of Proposition $4.8$. The size of the matrix $M^{-1}$ is irrevelant to their computation, and their formula can be generalized to our case to give the formula \begin{equation}
  	\hat{\z}_0(q)=(-1)^nq^{-\Delta} \sum_{r=0}^{\infty} \sum_{\epsilon \in  \{\pm 1\}^{n}} a_r \left( \prod_{j=1}^n \epsilon_j \right) q^{\frac{P}{4}\left(2r+(n-2)+\sum_{j=1}^n  \epsilon_j \frac{1}{p_j}\right)^2.}
  	\end{equation}
  	We now compute $\Psi(q).$  For $\lvert z \rvert <1$ we have
  	\begin{align}
  	G(z)= &\prod_{j=1}^{n-2} \left(z^{\frac{P}{p_j}}-z^{-\frac{P}{p_j}}\right)\left(z^{P}-z^{-P}\right)^{2-n}
  	\\ =& \sum_{r=0}^{\infty} \sum_{\epsilon \in  \{\pm 1\}^{n}} a_r \left( \prod_{j=1}^n \epsilon_j   \right) z^{2Pr+P(n-2)+\sum_{j=1}^{n}\epsilon_j  \frac{P}{p_j}}.
  	\end{align}
  	It follows that 
  	\begin{equation}
  	\Psi(q)=(-1)^n \sum_{r=0}^{\infty} \sum_{\epsilon \in  \{\pm 1\}^{n}} a_r  \left( \prod_{j=1}^n \epsilon_j  \right) q^{\frac{\left(2Pr+P(n-2)+\sum_{j=1}^{n}\epsilon_j  \frac{P}{p_j}\right)^2}{4P}}.
  	\end{equation}
  	This shows  \eqref{eq:numerodos}. \end{proof}
  
  We obtain the following corollary.
  \begin{corollary}
  	Let $\z_0 \in x^{-1/2}\C[[x^{-1}]]$ be the normalization of the Ohtsuki series from Theorem \ref{2.1}. We have an asymptotic expansion
  	\begin{equation}
  	q^{\Delta} \hat{\z}_0(X;q) \sim_{q\rightarrow 1} \frac{2\lambda}{\sqrt{\tau}} \z_0(1/\tau).
  	\end{equation}
  \end{corollary}

\begin{proof}
	This is a consequence of the integral formula \eqref{eq:zedhat} from Theorem \ref{thm:zedhat}
	\begin{align} q^{\Delta}\hat{\z}_0(X;q)&= \frac{\lambda}{\sqrt{\tau}} \int_{\varGamma} \exp(-\zeta / \tau )  \mathcal{B}(\z_0) (\zeta) \ \dv \zeta
	\\ &= \frac{\lambda}{\sqrt{\tau}} \sum_{\epsilon\in \{\pm1\}} \mathcal{L}_{\varGamma_{\epsilon}} (\mathcal{B}(\z_0))(1/\tau)
	\end{align}
	and Borel-Laplace resummation, which is stated as Theorem \ref{thm:BorelLaplaceresummation} below.  
	\end{proof}

  Let us now recall previous work on the $q$-series $\Psi.$ We start with the case  $n=3,$ for which more is known. As already mentioned in the introduction, Lawrence and Zagier have shown in \cite{LawrenceZagier} that the quantum invariant $\tau_k(X)$ can be recovered as the radial limit of $\Psi(q),$ as $q$ tends to $\exp(2\pi i/k)$.
   This was generalized to $n=4$ by Hikami in \cite{Hikami05b} but with corrections terms appearing.  The series $\Psi$ have interesting arithmetic properties; the coefficients $\chi(m)$ are periodic functions of period $2P$ and $\Psi$ is the so-called Eichler integral of a mock modular form with weight $3/2.$ As mentioned in the introduction the connection between quantum invariants and number theory was further pursued by Hikami in a number of articles \cite{Hikami05c, Hikami2004,  Hikami05b,Hikami06q, Hikami11}. For general $n\geq 3$ we mention again the work \cite{Hiroyuki} of Fuji, Iwaki, Murakami and Terashima which was discussed in the introduction.

  Let us now discuss what was previosly known about the $q$-series invariant $\hat{\z}_0(X).$ In \cite{GukovCiprian19} it was shown that when $X$ is a Briskekorn sphere $\Sigma(p_1,p_2,p_3),$ (i.e. $n=3$) then  $\hat{\z}_0(X)$ is a linear combination of so-called false theta functions.  The $q$-series invariant $\hat{\z}_0$ was also considered for certain Seifert fibered manifolds (with up to $n=4$ singular fibers) in the work \cite{ChungHJ18}, as well as a proposed analog of $\hat{\z}_0$ for higher rank gauge group - see also \cite{2019arXiv190913002P} for further developments in this direction. In this paper we work exclusively with $G=\SU(2)$).
  
  In connection with the work \cite{LawrenceZagier}, Zagier invented the notion of a quantum modular form. This notion was generalized by Bringmann et al. in \cite{BKM17}, where they introduce the notion of a higher depth quantum modular form.  For any $n\geq 3,$ it is known, that $\Psi$ is a linear combination of derivatives of quantum modular forms \cite{BMA18,BMA19june}. It is interesting to observe that $\Psi$ is obtained  from the Borel transform through a resummation process reminiscent of the median resummation of \cite{CostinGaroufalidis11}. Moreover as explained in \cite{ChengChunFerrariGukovHarrison18} it is expected that for a general $3$-manifold $M,$ Mock/false modular form duality is related to $\widehat{Z}_a(M;q)$, i.e. there exists an associated pair of a so-called Mock modular form and a so-called false modular form, and these are related by a $q \mapsto q^{-1}$ transformation and have the same transseries expression near $q \rightarrow 1.$ This is quite possibly connected to the conjecture $2$ in \cite{Garoufalidis08} (called the symmetry conjecture).  Let us also mention the work \cite{TudorGaroufalidis15} by Dimofte-Garoufalidis which connects modularity in quantum topology with complex Chern-Simons theory.
  
  \section{The asymptotic expansion of the GPPV invariant} \label{sec:radiallimit}
  
\begin{figure}[htp]
	\centering
	\includegraphics[scale=0.4]{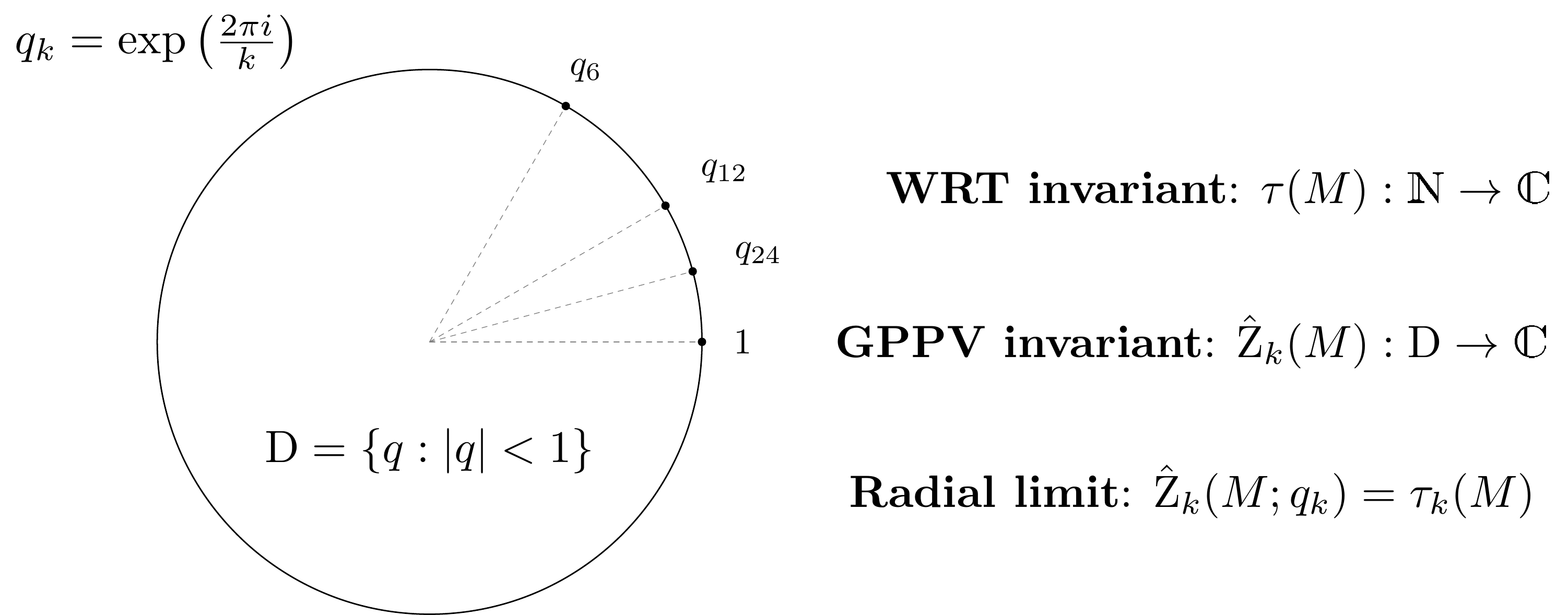}
	\caption{Analytic extension of $\tau_k(M)$.}
	\label{fig:analytic}
\end{figure}

The invention of $\hat{\z}$ was party motivated by an attempt to generalize the following discovery of Lawrence and Zagier. Set $q_k=\exp(2\pi i/k).$ For $n=3$ they proved in \cite{LawrenceZagier} the identity (for some $\sigma \in \Q)$  $$\tau_k\left(\Sigma(p_1,p_2,p_3)\right)(q_k-1)q_k^{\sigma}= -\frac{1}{2} \lim_{q \rightarrow q_k} \sum_{m=m_0}^{\infty} \chi_mq^{\frac{m^2}{4P}}.$$
For a closed oriented $3$-manifold $Y$ consider the normalized WRT invariant 
\begin{equation}
\z_{\CS}(Y;k)= \frac{\tau_k(Y)}{\tau_k(\Sph^2 \times \Sph^1)}.
\end{equation} 
We now state the radial limit conjecture.
\begin{conjecture}[\cite{GukovCiprian19}]  \label{Conj:radiallimt} Let $Y$ be a closed oriented $3$-manifold with $b_1(Y)=0.$ Set $$T=\Spin(Y)/\Z_2.$$ For every $a\in T,$ there exists invariants
	$$\Delta_a \in \Q, \ c \in \Z_+, \ \hat{\z}_a(q) \in 2^{-c}q^{\Delta_a}\Z[[q]],$$ with the following properties. The series $\hat{\z}_a(q)$ is convergent inside the unit disc $\{q\mid \  \lvert q \rvert <1\},$ and for infinitely many $k \in \N$ the radial limits $\lim_{q \uparrow \exp(2\pi i/k)} \hat{\z}_a(q)$ exists and we have that
	\begin{equation}
	\z_{\CS}(Y;k)=(i\sqrt{2k})^{-1} \sum_{a,b\in T} e^{2\pi i k lk(a,a)} \lvert \mathcal{W}_{b} \rvert^{-1} S_{a,b} \lim_{q\uparrow \exp(2\pi i/k)} \hat{\z}_b(q).
	\end{equation}
	Here $$S_{a,b}=\frac{e^{2\pi i kl(a,b)}+e^{-2\pi i kl(a,b)}}{\lvert \mathcal{W}_{a} \rvert \sqrt{\lvert \h_1(Y;Z)\rvert}},$$
	and $\mathcal{W}_x$ is the $\Z_2$-stabilizer of $x.$
\end{conjecture}
\begin{remark} Conjecture \ref{Conj:radiallimt} appeared in slighly different form in \cite{ChengChunFerrariGukovHarrison18, GukovPeiPutrovVada18,GukovMarinoPutrov}.
	\end{remark}  The level $k$ WRT invariant $\tau_k(M)$ of a closed oriented $3$-manifold $M$ can be seen as a function of the $k$-root of unity  $q_k=\exp(2\pi i/k)$, and as such it is a function of a certain subset of the boundary of the unit disc $\Disc=\{q : \lvert q \rvert <1\}.$ Assume $b_1(M)=0$ and define the $k$-dependent $q$-series
\begin{equation}
\hat{\z}_k(M;q)=(i\sqrt{2k})^{-1} \sum_{a,b\in T} e^{2\pi i k\  lk(a,a)} \lvert \mathcal{W}_{b} \rvert^{-1} S_{a,b}\hat{\z}_b(q).
\end{equation}
Then $\hat{\z}_k(M;q)$ is convergent for $q \in \Disc$ and the radial limit conjecture states
$$\lim_{q \uparrow q_k} \hat{\z}_k(M;q)=\tau_k(M).$$ Thus $\hat{\z}_k(M;q)$ can be seen as an analytic extension of $\tau_k(M)$ to the interior of the unit disc as illustrated in Figure \ref{fig:analytic} below

\subsection{Proof  of Theorem \ref{thm:radiallimit}}

To simplify notation, we write $$\hat{\z}_0(X;q)=\hat{\z}(q).$$ 
Recall the decomposition \eqref{phasedecomp} of the normalized quantum invariant $\tilde{\z}_k(X)$ into an integral part $\z^{\I}$ and a residue part $\z^{R}$.	In Lemma \ref{Lemmaet} we prove the existence of an analogous decomposition for $\hat{\z}_0(q)$ into a Laplace integral part $\mathcal{L}$ and a residue part $R$ 	\begin{equation} \label{eq:sumoftwoparts}
	q^{\Delta}\hat{\z}_0(q)= \mathcal{L}(\tau)+R(\tau),
	\end{equation}
	where we recall that $q = e^{2\pi i \tau}.$
We present in Proposition \ref{pro:ZagiersProposition} a standard result in complex analysis \cite{LawrenceZagier} which asserts that a $q$-series with periodic coefficients of mean value zero has an asymptotic expansion, as $q$ tends to a root of unity. We then show in Proposition \ref{pro:Fmeanvaluezero} that $R$ satisfy this hypothesis.  Finally, we apply Proposition \ref{pro:ZagiersProposition} to prove Theorem \ref{thm:radiallimit}.

\subsubsection{The decomposition of the GPPV invariant} 

Recall that $q=\exp(2 \pi i \tau)$ with $\tau \in \mathfrak{h}$ where $\mathfrak{h}$ denotes the upper half plane, and recall the definitions of $F \in \mathcal{M}(\C)$ and $g \in \C[x]$ given in \eqref{def:F}. Let $\tilde\Gamma_+ = e^{\frac{\pi i}{4}}$.
\begin{lemma} \label{Lemmaet}
Introduce the holomorphic functions $\mathcal{L},R \in \mathcal{O}(\mathfrak{h})$ given by
	\begin{align}
	\mathcal{L}(\tau) &= \frac{2\lambda}{\sqrt{\tau}} \int_{\varGamma_+} e^{-x /\tau} \mathcal{B}(\z_0)(x) \ \dv x,
	\\  R(\tau) &= -\frac{2 \lambda} {\sqrt{\tau} } \sum_{m=0 }^{\infty}  \Res\left( \exp\left( \frac{g(\xi)}{\tau}\right) F(\xi) ,\xi= 2 \pi i m\right).
	\end{align}
	Then we have that
	\begin{equation} \label{eq:sumoftwoparts}
	q^{\Delta}\hat{\z}_0(q)= \mathcal{L}(\tau)+R(\tau).
	\end{equation}
	For $\tau$   in the first quadrant $ \mathfrak{h}_+ = \{z\in \mathfrak{h} \mid \Re(z) >0\}$ we have that
	\begin{equation} \label{eq:DEDEFORMATION}
		\mathcal{L}(\tau) = \frac{2\lambda}{\sqrt{\tau}} \int_{0}^{\infty} e^{-x /\tau} \mathcal{B}(\z_0)(x) \ \dv x.
	\end{equation}
	
\end{lemma}

\begin{proof}[Proof of Lemma \ref{Lemmaet}]
	
 Recall the contour formula from Theorem \ref{thm:zedhat}

\begin{equation}
q^{\Delta} \hat{\z}_0(q) =\I(\tau) = \frac{\lambda}{\sqrt{\tau}} \int_{\varGamma} \exp(-x /\tau) \mathcal{B}(\z_0)(x) \ \dv x.
\end{equation}
Under the coordinate change
\begin{equation}
x=y^2
\end{equation} the contour $\varGamma$ corresponds to the contour
\begin{equation} 
\varPsi= e^{i\pi \frac38}\R_{+}+e^{i\pi \frac18}\R_+
\end{equation} 
Therefore we have that
\begin{equation}
q^{\Delta} \hat{\z}_0(q)= \frac{\lambda}{\sqrt{\tau}} \int_{\varPsi} \exp(-y^2 /\tau) \mathcal{B}(\z_0)(y^2)2y \ \dv y.
\end{equation} 
Introduce the meromorphic function $B \in \mathcal{M}(\C) $ given for all $y \in \C$ by \begin{equation} \label{eq:defB} B(y)=2y\mathcal{B}(\z_0)(y^2).
\end{equation}  By Theorem \ref{2.1} we have that
\begin{equation} \label{eq:B}
B(y)= \frac{2\kappa}{\pi i} G\left( \exp\left(\frac{\kappa y}{P}\right)\right)=\frac{2\kappa}{\pi i} F\left(\frac{\kappa y}{2}\right).
\end{equation}
From \eqref{eq:B} we see that $B$ is periodic with period $\kappa$, i.e. for all $m\in \Z$ we have that 
\begin{equation} \label{eq:periodicity}
B(y+\kappa m)=B(y).
\end{equation}

Let $\mathcal{P}$ be the set of poles of $B$. It follows from Theorem \ref{2.1} that $\mathcal{P}$ is a subset of the axis $e^{\pi i/4}\R$ and that 
\begin{equation}
 \{ -\omega^2  \mid  \omega \in \mathcal{P}  \}  = 2\pi i\CS_{\C}(X) + 2\pi i\Z.
\end{equation}  
Write $$\varPsi_\pm= e^{i \pi (\frac14 \mp \frac18)}\R_+.$$ We will now  apply Cauchy's residue theorem to move $\varPsi_-$ across $e^{i \pi /4}\R_+$ to $\varPsi_+$  in order to obtain the formula \eqref{eq:sumoftwoparts}. Deform $\varPsi_\pm$ on the complement of a neighbourhood around the origin to two curves $L_\pm$, which are parallel to $e^{i \pi /4}\R_+$ outside this neighbourhood of the origin, as indicated in Figure \ref{fig:ContourL}. Set $$L=L_+ \cup L_-.$$ We first show that 
\begin{equation}  
\label{eq:shiftofContour} \int_{\varPsi} \exp(-y^2/\tau) B(y) \ \dv y= \int_{L} \exp(-y^2/\tau) B(y) \ \dv y,   \end{equation}
and then we show that the right hand side of \eqref{eq:shiftofContour} can be rewritten as a sum of residues. Let $R>0$ be a positive constant, and let $R_\pm$ be the arc segment of the circle of radius $R$, which connects $\varPsi_\pm$ and $L_\pm$.  Because $L_\pm$ is parallel to $\kappa \R_+$ outside a neighbourhood of the origin, there exists a real positive constant $b_0 >0$ such that every $y \in R_\pm$ is of the form
$$y= \kappa a \pm b, (a,b) \in \R_+ \times [b_0,+\infty)$$
and therefore exists a positive real constant $A>0$ independent of $R$, which gives an upper bound
\begin{equation}
\left\lvert \frac{1}{1-e^{-\kappa y}}\right\rvert <\frac{1}{1-\exp\left(-\Re(\kappa) b_0\right)}:=A>0
\end{equation}
for all $y \in R_+$.
It follows that we obtain a uniform estimate 
\begin{align} \label{eq:Estimation1} \begin{split} B(y) &= 2^{n-2}\frac{2\kappa}{\pi i}\frac{\prod_{j=1}^n \sinh(\frac{\kappa y}{2p_j})}{\left( e^{\kappa y/2}(1-e^{-y\kappa }) \right)^{n-2}}
\\&=  O \left( e^{-y\kappa /2}\prod_{j=1}^n \sinh\left(\frac{\kappa y}{2p_j}\right) \right)= O(e^{A_1R}). \end{split}  \end{align}
for all $y \in R_+$ for a real constant $A_1$.
For fixed  $\tau$, there exists a positive real constant $A_2>0$ giving a uniform bound on $R_+$
\begin{equation}  \label{eq:Estimation2}
e^{-y^2/\tau}=O(e^{-A_2R^2}).
\end{equation}
By combining the estimates \eqref{eq:Estimation1} and \eqref{eq:Estimation2}, and using that the arc length of $R_\pm$ is proportional to $R$, we obtain the estimate
\begin{equation}
\int_{R_+} \exp(-y^2/\tau) B(y) \ \dv y= O(Re^{-A_2R^2+A_1R}).
\end{equation}
By similar reasoning, there exists constants $B_2>0,B_1 \in \R$ giving the estimate
\begin{equation}
\int_{R_-} \exp(-y^2/\tau) B(y) \ \dv y= O(Re^{-B_2R^2+B_1R}).
\end{equation}
Thus we obtain that
\begin{equation} \label{eq:LimitZero}
\lim_{R \rightarrow \infty} \sum_{\epsilon \in \{\pm1\}} \int_{R_\epsilon} \exp(-y^2/\tau) B(y) \ \dv y=0,
\end{equation}
which gives the desired identity \eqref{eq:shiftofContour} .

 We now turn to the computation of $\int_{L}  \exp(-y^2/\tau) B(y) \ \dv y$. For each $m\in \N$, let $L_m$ be a small line segment with 
 $$L_m \cap e^{i\pi/4}\R_+=\{m \kappa\},$$ 
 and which meets $e^{i\pi/4}\R_+$ in a right angle. We can arrange that $L_m$ is of fixed lenght and that $L_m $ meet $L_\pm$ in a point. Thus we have
\begin{equation} \label{eq:periodicity2} L_m=L_0+m\kappa
\end{equation} Let $L^m_\pm$ be the bounded component of $L_\pm \setminus L_m,$ and let $\mathcal{P}_m \subset \mathcal{P}$ be the set of poles of $B$ that lie within the bounded component of the complement of the contour $$\varPsi_m = L^m_+ \cup L_m  \cup L^m_-.$$ See Figure \ref{fig:ContourL}.
\begin{figure}[htp]
	\centering
	\includegraphics[scale=0.9]{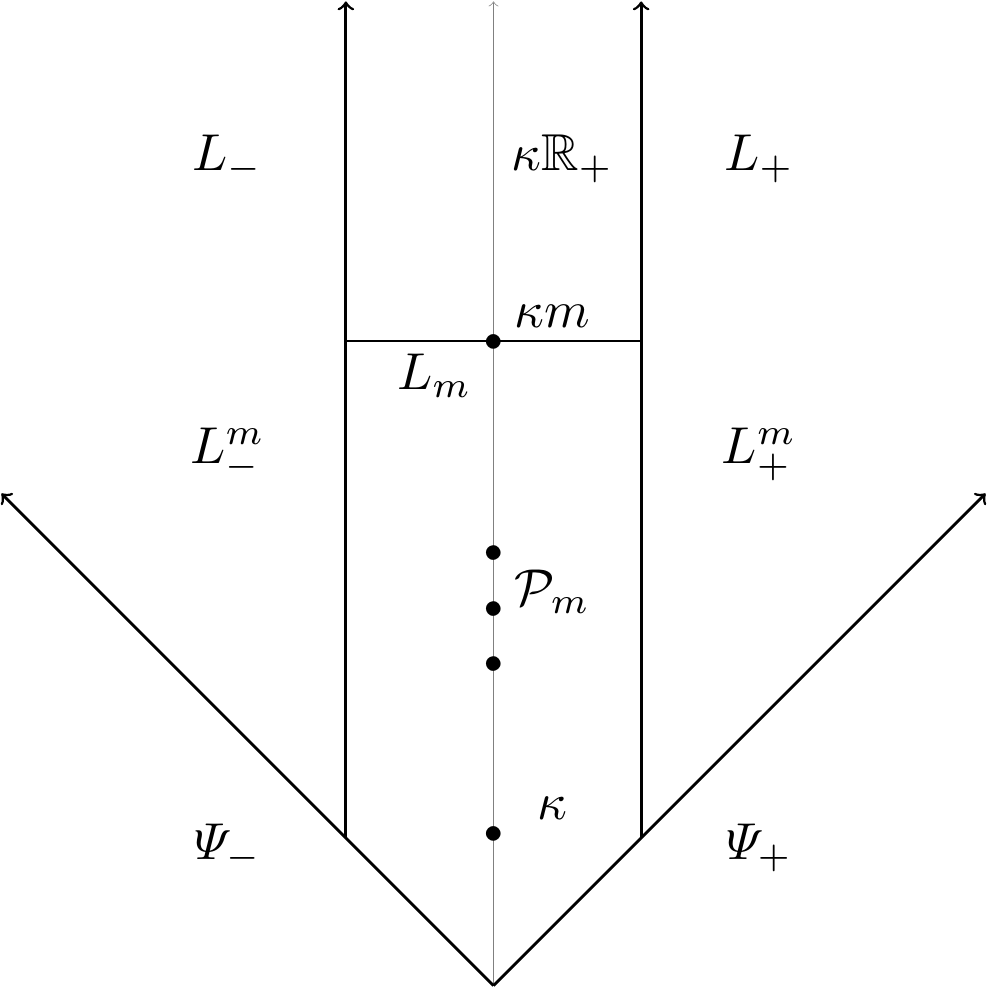}
	\caption{The contours $L_\pm, L^m_\pm$ and $ \varPsi_\pm$ and the subset of poles $\mathcal{P}_m$ drawn rotated, such that $\kappa \mathbb{R}_+$ points straight up.}
	\label{fig:ContourL}
\end{figure}   Equip $\varPsi_m$ with the counter clockwise orientation. An application of Cauchy's residue theorem now gives
\begin{multline} \label{eq:dabel} \begin{split}
2\pi i\sum_{\omega \in \mathcal{P}_m} \Res(e^{-y^2/\tau}B(y), y=\omega)  =  \int_{\varPsi_m} \exp(-y^2/\tau) B(y) \ \dv y
\\ =\int_{L_m} \exp(-y^2/\tau) B(y) \ \dv y +\sum_{\epsilon \in \{\pm1\}} \epsilon  \int_{L^m_\epsilon} \exp(-y^2/\tau) B(y) \ \dv y. \end{split} 
\end{multline}
Because $B$ is $\kappa$ periodic as stated in \eqref{eq:periodicity} and $L_{m}=L_0+ m\kappa$, we see that that there exists $C>0$ giving a uniform bound
\begin{equation}
C> \sup \{ \lvert B(y) \rvert \ \mid \ y \in \cup_{m\in \N} L_m \}.
\end{equation}
Because of this universal bound, it is easy too see that
\begin{equation}
\lim_{m \rightarrow \infty} \int_{L_m} \exp(-y^2/\tau) B(y) \ \dv y=0.
\end{equation}
It follows that the right hand side of \eqref{eq:dabel} converges to 
$$\sum_{\epsilon \in \{\pm1\}} \int_{L_\epsilon} \exp(-y^2/ \tau) B(y)  \ \dv y=\sum_{\epsilon \in \{\pm1\}} \epsilon  \int_{\varPsi_\epsilon} \exp(-y^2/ \tau) B(y)  \ \dv y.$$ 
This also implies that the sum of residues is convergent. 

Let us now recall a simple transformation law for residues. Let $z_0\in \C$ and let $f\in \mathcal{M}_{z_0}(\C)$ be the germ of a meromorphic function with a pole at $z_0.$ Assume $w_0 \in \C$ and that $z \in \mathcal{O}_{w_0}(\C)$ satisfies $z(w_0)=z_0, \dot{z}(w_0)\not=0.$ If either $z_0$ is a simple pole, or $z(w)$ is linear in $w,$ then we have that
\begin{equation} \label{ResiduePropositionFormula}
\text{Res}(f(z(w)), w=w_0)= \frac{\text{Res}(f(z),z=z(w_0))}{\dot{z}(w_0)}.
\end{equation} 
 Introduce the variable
\begin{equation} \label{eq:xi}
\xi= \frac{\kappa y}{2}.
\end{equation} By using the relation \eqref{eq:B} between $B$ and $F$, the relation \eqref{eq:defB}  between $B$ and $\mathcal{B}(\z_0)$ and the tranformation law \eqref{ResiduePropositionFormula} we obtain
\begin{align} \begin{split} 
q^{\Delta} \hat{\z}_0(q)& = \frac{\lambda}{\tau} \int_{\varPsi} \exp\left(-\frac{y^2}{\tau}\right) B(y) \ \dv y
\\&=\frac{2\lambda}{\tau} \left( \int_{\varPsi_+} \exp\left(-\frac{y^2}{\tau}\right) B(y) \ \dv y -  2\pi i\sum_{\omega \in \mathcal{P}} \Res(e^{-y^2/\tau}B(y), y=\omega)  \right)
\\ &= \frac{2\lambda}{\tau} \left( \int_{\Gamma_+} \exp\left(-\frac{x}{\tau}\right) \mathcal{B}(\z_0)(x) \ \dv x -  \sum_{m=1}^{\infty} \Res(e^{-g(\xi)/\tau}F(\xi), \xi=2\pi i m) \right)
\\ &= \mathcal{L}(\tau)+R(\tau). \end{split}
\end{align}
Finally, we prove \eqref{eq:DEDEFORMATION} for $\tau \in \mathfrak{h}_+$. First observe that for $x$ and $\tau$ in the upper right half plane we have
\begin{equation} \label{Positve}
\Re(x/\tau) >0.
\end{equation} Push the contour $\tilde{\Gamma}_+$  to $\R_+$. If the integral is invariant under this deformation of the contour, we obtain the desired identity. To see that the  integral is invariant under this deformation of the contour, we apply a limiting argument, together with Cauchy's residue formula. To that end, let $R>0$ be a positive parameter, and let $C_R$ be the  arc segment of the circle of radius $R$, which connects $R$ to $Re^{i\pi \frac{3}{4}}$ and stays in the upper half plane. As we are not moving the contour across any singularities of $\mathcal{B}(\z_0),$ the only difficulty is to show that 
\begin{equation} \label{eq:LIIMITzero}
\lim_{R\rightarrow +\infty} \int_{ C_R} e^{-x / \tau} \mathcal{B}(\z_0)(x) \ \dv x=0.
\end{equation}
As $C_R$ remain at least a fixed distance away from the axis of poles of $\mathcal{B}(\z_0)(x)$, the limit \eqref{eq:LIIMITzero} follows by \eqref{Positve} together with arguments similar to the arguments giving the limit \eqref{eq:LimitZero} above.  \end{proof}

\subsubsection{Asymptotic expansions of $q$-series with periodic coefficients}

Let $B_m(x)$ denote the $m$-th Bernoulli polynomial, i.e.
\begin{equation} \label{def:Bernoulli}
\frac{te^{tx}}{e^t-1}= \sum_{m=0}^{\infty} \frac{B_m(x)}{m!} t^m.
\end{equation} We recall the following result.

\begin{proposition}[\cite{hikami2003qseries,LawrenceZagier}] \label{pro:ZagiersProposition}
	Let $C: \Z \rightarrow \C$ be a periodic function with period $M$ and mean value equal to zero
	\begin{equation}
	\sum_{m=1}^M C(n)=0.
	\end{equation}
	Consider the $L$-series $L(s,C)$, which for $\Re(s)>1$ is defined by 
	\begin{equation} \label{eq:Lseries}
	L(s,C)= \sum_{m=1}^{\infty} \frac{C(m)}{m^s}.
	\end{equation} 
	This $L$-series admits an analytic extension to all of $\C$ and for $r \in \N$ 
	\begin{equation} \label{eq:Lseriesr}
	L(-r,C)= -\frac{M^r}{r+1} \sum_{m=1}^M C(m) B_{r+1}\left(\frac{m}{M}\right).
	\end{equation}
	For any polynomial $Q$ of degree $d$
	$$Q(x)=\sum_{u=0}^d q_{u} x^u \in \C[x]$$ the following asymptotic expansions hold for real and positive $t$
\end{proposition}
	
	\begin{equation} \label{eq:qseriesexpansion}
	\sum_{m=1}^{\infty} e^{-tm^2} C(m) Q(m) \ \underset{t \rightarrow 0}{\sim}  \sum_{u=0}^d \sum_{r=0}^{\infty} q_u L(-2r-u,C) \frac{(-t)^r}{r!}.
 	\end{equation}
	
	\begin{proof} 
		The existence of the analytic extension of the $L$-series of $C$, as well as the explicit evaluation \eqref{eq:Lseriesr} are proven in \cite{LawrenceZagier}.
		
		 In \cite{hikami2003qseries,LawrenceZagier} the following asymptotic expansions are proven\begin{align} \begin{split} \label{eq:asymptotixexpansions} 
		 &	\sum_{m=1}^{\infty} C(m) e^{-tm^2} \ \underset{t \rightarrow 0}{\sim} \  \sum_{r=0}^{\infty} L(-2r,C) \frac{(-t)^r}{r!},\\   &\sum_{m=1}^{\infty} m C(m) e^{-tm^2} \ \underset{t \rightarrow 0}{\sim} \  \sum_{r=0}^{\infty} L(-2r-1,C) \frac{(-t)^r}{r!}. \end{split} 	\end{align}	We have	\begin{equation} \label{eq:difff}	\sum_{m=1}^{\infty} e^{-tm^2}C(m)Q(m)= \sum_{j=0}^{\infty} \frac{\partial^j}{\partial (-t)^j} \sum_{m=0}^{\infty} (q_{2j+1}m+q_{2j})e^{-tm^2},	\end{equation}
		 where it is understood that $q_l=0$ for $l>d$.	The expansion \eqref{eq:qseriesexpansion} follows formally from differentiating the expansions given in \eqref{eq:asymptotixexpansions}. This differentiation is valid because Poincare asymptotic expansions of analytic functions which are valid on suitable sectors can be termwise differentiated. Clearly $t \mapsto \sum_{m\geq 0} C(m)Q(m)\exp(-tm^2)$ is an analytic function of $t$ in a small tubular neighbourhood of $(0,1]$, and from the proof given in \cite{LawrenceZagier} it is clear that the asymptotic expansions \eqref{eq:asymptotixexpansions} are valid on such a small sector. 
	\end{proof}

Recall the definition of the meromorphic function $F$ given in  \eqref{def:F}. Next we prove that the coefficients of the principal part of $F$ at poles are periodic functions with mean value equal to zero.
\begin{proposition} \label{pro:Fmeanvaluezero}
	For $j=1,2,...,n-2$ define $f_j: \Z\rightarrow \C$ as the coeficients of the principal part of $F$ at $2\pi i m$ for $m \in \Z$, e.g. for $y$ near $2\pi im$
	\begin{equation}
	F(y) = \sum_{j=1}^{n-2} f_j(m) (y-2\pi i m) ^{-j}+\rm{reg.}.
	\end{equation}
	Then each $f_j$ is $2P$-periodic and if $P$ is even, then we have for each even $k \in \Z$ 
	\begin{equation} \label{eq:MeanValueZero}
	\sum_{m =1,...,2P } e^{kg(2\pi im)} f_j(m)=0.
	\end{equation}
	\end{proposition}

\begin{proof}
The periodicity of the functions $f_j,j=1,...,n-3$ follow directly from the $4\pi i P$-periodicity of $F$. 

We now prove  \eqref{eq:MeanValueZero} assuming $P$ is even - which is equivalent to exactly one the $p_j$ being even. Using the definition \eqref{def:F} of $F$ we obtain
\begin{align}
F(y + 2\pi iP)&=\frac{1}{4} \sinh(y/2 + \pi i P)^{2-n} \prod_{j=1}^{n} \sinh(y/p_j +\pi i P/p_j) 
\\ &= (-1)^{nP} (-1)^{\sum_{j=1}^n \frac{P}{p_j}} F(y)=-F(y).
\end{align}
This implies that for each $j=1,...,n-3$ and $m=1,...,P$ we have
\begin{equation} f_j(m+P)=-f_j(m). \end{equation}
On the other hand we have for integral $m$
\begin{align}
g(2\pi i(P + m))&=\frac{i(2\pi i(P+ m))^2}{8\pi P}= \frac{-i\pi m^2}{2P}\mp i\pi mP-\frac{i\pi P}{2} 
\\ &=g(2\pi i(P\mp m)) \mod \pi i \Z.
\end{align} It follows that for even $k$ we have a pairwise cancellation
\begin{align}
&\sum_{m=1}^{2P} e^{kg(2\pi i m)} f_j(m)= \sum_{m=1}^{P} e^{kg(2\pi i m)} f_j(m)+e^{kg(2\pi i (m+P))} f_j(m+P)
\\&=\sum_{m=1}^{P} e^{kg(2\pi i m)} f_j(m)- e^{kg(2\pi i m)} f_j(m)=0.
\end{align} This concludes the proof. \end{proof}

\subsubsection{The asymptotic expansion of the GPPV invariant}

Recall the definition the functions $F$ and $g$ given in \eqref{def:F}. Let $m\in \N$. For $y$ close to $2 \pi im$ we use the notation of Propositon \ref{pro:Fmeanvaluezero} and write
$$F(y)= \sum_{j=1}^{n-2} f_j(m) (y-2\pi im)^{-j}+\rm{reg.},$$
where each $f_j: \Z\rightarrow \C$ is $2P$-periodic. For each $l\in \N$ there exists a uniquely determined polynomial \begin{equation} \label{def:Taylorcoefficients}
P_l(x,y) \in Q[\pi i,x,y]
\end{equation} such that 
\begin{equation} \label{eq:Taylorcoefficients}
\frac{1}{l!}\frac{\partial^l e^{g(y)/\tau}}{\partial y^l}(2\pi i m)= e^{g(2\pi im)/\tau} P_l(\tau^{-1},m).
\end{equation}
There exists uniquely determined complex coefficients $p_{l,u,v,}$ with
\begin{equation}
P_{l}(\tau^{-1},m)= \sum_{u,v} p_{l,u,v} \tau^{-u} m^v.
\end{equation}
Recall that $B_m(x)$ denotes the $m$'th Bernoulli polynomial and is defined by \eqref{def:Bernoulli}. Define for each $\theta \in \CS_{\C}(X)$ the following  polynomials with coefficients in power series
\begin{align} \begin{split}  \label{def:Rtheta}
R_{\theta}(k,t) &= \sum_{\substack{m=1,...,2P, \\ g(2\pi im) \equiv 2\pi i\theta }} \sum_{j,u,v}  f_j(m) k^up_{j-1,u,v}  \sum_{r=0}^{\infty}  \frac{(2P)^{2r+v}}{2r+v+1}  B_{2r+v+1}\left(\frac{m}{2P}\right) \frac{(-t)^r}{r!},
\\ \check{\z}_{\theta}(k,t) &= \sum_{\substack{m=1,...,2P, \\ g(2\pi im) \equiv 2\pi i\theta }} \sum_{j,u,v}  f_j(m) k^up_{j-1,u,v}  \sum_{r=1}^{\infty}  \frac{(2P)^{2r+v}}{2r+v+1}  B_{2r+v+1}\left(\frac{m}{2P}\right) \frac{(-t)^r}{r!}
\end{split} 
\end{align}
Observe $$\check{\z}_{\theta}(k,t)=R_{\theta}(k,t)-R_{\theta}(k,0),$$
i.e. $\check{\z}_{\theta}$ is equal to $R_{\theta}$ minus the constant in the parameter $t$.

We can now prove Theorem \ref{thm:radiallimit}.

\begin{proof}[Proof of Theorem \ref{thm:radiallimit}]
	
 Recall the decompsition
	$$\frac{\sqrt{\tau}}{\lambda}q^{\Delta}\hat{\z}=\mathcal{L}(\tau)+R(\tau)$$ given in Equation \eqref{eq:sumoftwoparts} in Lemma \ref{Lemmaet}. Recall the decomposition of the normalized quantum invariant
	$$\z_k(X)= \z^{\I}(k)+\z^R(k)$$ given in \eqref{phasedecomp}. This decomposition together with equation \eqref{def:normquant}, which relates the normalized quantum invariant $\tilde{\z}_k(X)$ to the WRT invariant $\tau_k(X)$, shows the radial limit idenity can be proved by proving the following two limits. 
	\begin{align} & \z^{\I}(k) = \lim_{\tau \uparrow 1/k} \frac{\sqrt{\tau}}{2\lambda} \mathcal{L}(\tau),
 & \z^R(k) = \lim_{\tau \uparrow 1/k}, \frac{\sqrt{\tau}}{2\lambda}R(\tau).\end{align}  
 Observe that as $X$ is an integral homology sphere, the only $\Spin$ structure is $a=0$, and therefore the radial limit conjecture  reduces (up to a scalar factor) to equation \eqref{eq:radiallimit}.  We first focus on the integral part $\mathcal{L}(\tau)$. For every $k\in \N^*$ the integral part  $\mathcal{L}(\tau)$ extends continuosly to $\tau=1/k$ and it follows from Equations \eqref{phasedecomp}, \eqref{ComputeBorel} and \eqref{eq:DEDEFORMATION}  that 
		\begin{equation}
		(2\lambda \sqrt{k})^{-1} \mathcal{L}(1/k)= \mathcal{L}_{\R_+}(\mathcal{B}(\z_0))(k)=\z^{\I}(k).
		\end{equation}
		Now recalling that $\tau_{k,t}=(k-i\frac{2Pt}{\pi})^{-1}$, we see that the non-trivial parts left in order to prove  the asymptotic expansion \eqref{eq:ZedhatExpansion} stated in Theorem \ref{thm:radiallimit} are the asymptotic expansion 
		\begin{equation} \label{eq:MASTEREXPANSION}
		\frac{\sqrt{\tau_{k,t}}}{2\lambda} R(\tau_{k,t}) \ \underset{ t \rightarrow 0}{\sim} \ \sum_{\theta \in \CS_{\C}(X)} e^{2\pi ik\theta} R_{\theta}(k,t),
		\end{equation} 
		and the identity
		\begin{equation} \label{eq:FINALIDENTITY} \sum_{\theta \in \CS_{\C}(X)} e^{2\pi ik\theta} R_{\theta}(k,0)= \z^{R}(k).
		\end{equation}
We start with the expansion \eqref{eq:MASTEREXPANSION}. To ease notation we set
\begin{equation}
g(2\pi im)=\theta_m.
\end{equation} We have
\begin{equation}
\Res(F(y)e^{g(y)/\tau}, 2\pi im)=  e^{\theta_m/\tau} \sum_{j=1}^{n-3} f_j(m) P_{j-1}(\tau^{-1},m).
\end{equation}
and accordingly
\begin{equation}
\frac{\sqrt{\tau}}{2\lambda}R(\tau)=-\sum_{m =0}^\infty  e^{\theta_m/\tau} \sum_{j=1}^{n-3} f_j(m) P_{j-1}(\tau^{-1},m),
\end{equation}
Note that $\tau_{k,t} \in \mathfrak{h}$ for $t\in (0,1)$.
The function $g$ has the following transformation property\begin{equation}g(y+4\pi i m P)=\frac{i(y+4\pi i mP)^2}{8\pi P}=g(y)-my- 2 \pi i  m^2 P.\end{equation} Because of this, and the $2P$ periodicity of the $f_j$'s, we can write
\begin{align} \begin{split}  \label{MELLEMREGNING}
\frac{\sqrt{\tau}}{2\lambda}R(\tau)= &-\sum_{m=0}^{\infty} e^{- t m^2} e^{kg(2\pi i m)} \sum_{j=1}^{n-3} f_j(m) P_{j-1}(\tau^{-1},m)
\\ = &- \sum_{j=1}^{n-2} \sum_{u,v} \sum_{m=0}^{\infty} e^{- t  m^2} e^{k\theta_m}  f_j(m) p_{j-1,u,v} \tau^{-u}m^v.
\end{split} 
\end{align}
Now for each $j,u=0,...,n-3$ we can apply Proposition \ref{pro:ZagiersProposition} to the $2P$-periodic function of mean value zero given by
\begin{equation}
C_j(m)=e^{k\theta_m}  f_j(m)
\end{equation}
and the polynomial
\begin{equation}
Q_{j,u}(m)= \sum_v p_{j-1,u,v} m^v.
\end{equation}
The fact that each $C_j$ is of mean value equal to zero follows from Proposition \ref{pro:Fmeanvaluezero}. The result of applying Proposition \ref{pro:ZagiersProposition} is
\begin{align}
\frac{\sqrt{\tau}}{2\lambda}R(\tau) \ & \underset{t \rightarrow 0}{\sim } \ -\sum_{j,u,v} \sum_{r=0}^{\infty}  k^up_{j-1,u,v} L(-2r-v,C_j) \frac{(-t)^r}{r!}
\\ &= \sum_{j,u,v} \sum_{r=0}^{\infty}  k^up_{j-1,u,v} \frac{(2P)^{2r+v}}{2r+v+1} \sum_{m=1}^{2P} C_j(m) B_{2r+v+1}\left(\frac{m}{2P}\right) \frac{(-t)^r}{r!}
\\ &= \sum_{m=1}^{2P} \sum_{j,u,v} e^{k\theta_m}  f_j(m) k^up_{j-1,u,v}  \sum_{r=0}^{\infty}  \frac{(2P)^{2r+v}}{2r+v+1}  B_{2r+v+1}\left(\frac{m}{2P}\right) \frac{(-t)^r}{r!}.
\end{align}
This establishes the asymptotic expansion \eqref{eq:MASTEREXPANSION}.

We now turn to the identity \eqref{eq:FINALIDENTITY}.  Set
\begin{align} \begin{split}  \label{eq:MUAHAH}
R_0(k)&=\sum_{\theta \in \CS_{\C}(X)} e^{2\pi ik\theta} R_{\theta}(k,0)= \sum_{m=1}^{2P} \sum_{j,u,v} e^{k\theta_m}  f_j(m) p_{j-1,u,v} k^u    \frac{(2P)^{v}}{v+1}  B_{v+1}\left(\frac{m}{2P}\right).
\end{split}
\end{align}
Let $x$ be a complex coordinate near $0$ and set $y_m=x+2\pi i m$ for $m=1,2,...,2P$. Recall that
\begin{equation}
\z^R(k)=- \sum_{m=1}^{2P-1}  \Res\left(  \frac{F(y)\exp\left(kg(y)\right)}{1-\exp(-ky)} , y=2\pi im \right).
\end{equation} 
We have that
\begin{align} \label{eq:vigtigdel}
g(y_m)=g(x+ 2\pi im)=\frac{i(x+2\pi i m)^2}{8\pi P}= g(2\pi im)+g(x)-\frac{xm}{2P},
\end{align}
and accordingly
\begin{equation} \label{eq:VIGTIGDEL}
\exp(kg(y_m) )=\exp\left(k\theta_m+kg(x) \right) \cdot \exp\left(- k\frac{xm}{2P}\right).
\end{equation}
Recall that the polynomials $P_l$ defined in \eqref{eq:Taylorcoefficients} as the coefficients of the Taylor series of $e^{kg(y)}$. Therefore it follows from Cauchy's formula for multiplication of power series, the formula for the Taylor expansion of the exponential and the identity \eqref{eq:VIGTIGDEL} that the following  holds for all $k ,m$
\begin{equation} 
P_c(k,m)= \sum_{a+b=c} P_a(k,0)\left(\frac{-mk}{2P}\right)^b \frac{1}{b!}.
\end{equation}
Writing this out in terms of coefficients gives
\begin{align}
\sum_{c,u,v} p_{c,u,v} k^u m^v &= \sum_{a+b=c} \sum_{s} p_{a,s,0} \left( \frac{(-1)}{2P} \right)^b \frac{1}{b!} k^{s+b} m^b.
\end{align}
This is equivalent to the identities
\begin{equation} \label{eq:MAASTERIDENTITY}
p_{j,u,0} \left( \frac{-1}{2P} \right)^v \frac{1}{v!}= p_{j+v,u+v,v}.
\end{equation}
Recall the definition \eqref{def:Bernoulli} of the Bernoulli polynomials $B_m$. Write
\begin{align}
&- \frac{F(y_m)\exp\left(kg(y_m)\right)}{1-\exp(-ky_m)} =F(y_m) \exp\left(k(g(x)+\theta_m\right) \frac{\exp\left(\frac{-kxm}{2P} \right)}{\exp(-kx)-1}=
	\\ &\left(\sum_{j=1}^{n-3} f_j(m)x^{-j}+\rm{reg.}\right)\left(e^{k\theta_m}\sum^{\infty}_{a=0} P_a(k,0)x^a\right)\left(\sum_{b=0}^{\infty} \frac{B_b\left(\frac{m}{2P}\right)}{b!} (-kx)^{b-1} \right).
\end{align}
By comparing with Equation \eqref{eq:MUAHAH} and using the facts that $F$ has a multiple order zero at multiples of $4\pi iP$ and that we know that the $k^{-1}$ term cancels in $\z^R(k)$, we obtain the desired identity
\begin{align}
&\z^R(k)= \sum_{m=1}^{2P} \sum_{j=1}^{n-2} \sum_{\substack{a+b=j, \\ a\geq 0, b \geq 1}} e^{k\theta_m} f_j(m)P_a(k,0)\frac{B_b\left(\frac{m}{2P}\right)}{b!} (-k)^{b-1}
\\ &=  \sum_{m=1}^{2P} \sum_{j=1}^{n-2}  \sum_{\substack{a+b=j, \\ a\geq 0, b \geq 1}} e^{k\theta_m} f_j(m)P_a(k,0)\left(\frac{-1}{2P}\right)^{b-1} \frac{k^{b-1}}{(b-1)!}  \frac{\left(2P\right)^{b-1}}{b}B_b\left(\frac{m}{2P}\right)
\\ &= \sum_{m=1}^{2P} \sum_{j=1}^{n-2}  \sum_{\substack{a+b=j, \\ a\geq 0, b \geq 1}} \sum_{s} e^{k\theta_m} f_j(m)p_{a,s,0}\left(\frac{-1}{2P}\right)^{b-1} \frac{k^{s+b-1}}{(b-1)!}  \frac{\left(2P\right)^{b-1}}{b}B_b\left(\frac{m}{2P}\right)
\\ \begin{split}  \label{eq:UHAHAHA} &=\sum_{m=1}^{2P} \sum_{j=1}^{n-2}  \sum_{\substack{s,a+b=j, \\ a\geq 0, b \geq 1}} \sum_{s} e^{k\theta_m} f_j(m) p_{a+b-1,s+b-1,b-1}  k^{s+b-1}   \frac{(2P)^{b-1}}{b} B_b\left(\frac{m}{2P}\right) \end{split}
\\ &= \label{eq:UAA} \sum_{m=1}^{2P} \sum_{j=1}^{n-2}  \sum_{u,v}  e^{k\theta_m} f_j(m) p_{j-1,u,v} k^u \frac{(2P)^{v}}{v+1} B_{v+1}\left(\frac{m}{2P}\right)
\\&=R_0(k). 
\end{align}
In \eqref{eq:UHAHAHA} we used the identity \eqref{eq:MAASTERIDENTITY}, and in \eqref{eq:UAA} we set $u=s+b-1$ and $v=b-1$. This finishes the proof. \end{proof}

\begin{appendices}

\section{Resurgence and resummation}  \label{sec:Borel}

\subsection{Resurgent functions and the Borel transform}

The theory of resurgence was originally developed by Écalle in \cite{Ecalle81a} and \cite{Ecalle81b}.  See \cite{MitschiSauzin16} for an introduction to the mathematical theory of resurgence and see \cite{DunneUnsal15} for an introduction to the general use of resurgence in quantum field theory. Garoufalidis \cite{Garoufalidis08} and Witten \cite{Witten11} where the pioneers of the use of resurgence in quantum Chern-Simons theory. \begin{definition} \label{def:res} For a Riemann surface $C$ with universal covering space 
	$$ \pi: \tilde{C}\rightarrow C$$ the algebra of resurgent functions is  $\mathcal{R}(C)=\mathcal{M}(\tilde{C}).$
\end{definition}

One source of resurgent functions are the Borel transforms of Laplace integrals. We now introduce the Borel transform. Let $\Gamma \in \mathcal{M}(\C)$ be the Gamma function, which for $z \in \C$ with $\re(z)>0$ is defined by
\begin{equation}
\Gamma(z)= \int_0^{\infty} e^{-t} t^{z-1} \dv t.
\end{equation}

\begin{definition} \label{def:Bor}  Let $\{\alpha_j \}_{j=0}^{\infty}$ be an increasing sequence of positive real numbers,  $\{\beta_j \}_{j=0}^{\infty}$ a sequence of non-negative integers and $\{c_j\}_{j=0}^{\infty}$ a sequence of complex numbers. Consider  the formal series $$\tilde{\varphi}(\lambda)=\sum_{j=0}^{\infty} c_j \lambda^{-\alpha_j} \log(\lambda)^{\beta_j}.$$ The Borel transform of $\tilde{\varphi}(\lambda)$ is given by the formal series
	\begin{equation} \label{def:Borel}
	\mathcal{B}(\tilde{\varphi})(\zeta) =\sum_{j=0}^{\infty} c_j (-1)^{\beta_j} \frac{\partial^{\beta_j}}{\partial {\alpha_j}^{\beta_j}} \left( \frac{\zeta^{\alpha_j-1}}{\Gamma(\alpha_j)}  \right).
	\end{equation}
\end{definition}

\subsection{Borel-Laplace resummation}We now discuss in more detail the relation between the Borel transform $\mathcal{B}$ and the Laplace transform, which we now introduce. Let $\gamma \subset \C$ be an oriented contour. Let $g$ be a measurable function defined in a neighbourhood of $\gamma$. Denote by $\mathcal{L}_{\gamma}(g)$ the Laplace transform given by
\begin{equation} \label{eq:Laplacetransform}
\mathcal{L}_\gamma (g) (\lambda) = \int_\gamma \exp(- \lambda \cdot  z ) g(z) \ \dv z,
\end{equation}
for all $\zeta  \in \C$ such that the integral is absolutely convergent. Here  we think of $\lambda \in \C^*$ as a large modulus asymptotics parameter.
 For any $\alpha \in \C^*$ we let the contour $\alpha \R_+$ be oriented in the direction of $\alpha$ unless we state otherwise.

That the transforms $\mathcal{B}$ and $ \mathcal{L}_{\R_+}$ are formally inverses of each other should be understood as follows. If $\alpha \in \C$ satisfies $\re(\alpha)>-1$ and $l\in \N$ then
\begin{equation} \label{eq:Laplacedas}
\mathcal{L}_{\R_+}\left(\zeta^{\alpha}\log(\zeta)^l\right)= \frac{\text{d}^l}{\dv \alpha^l} \left(\frac{\Gamma(\alpha+1)}{\lambda^{\alpha+1}} \right).
\end{equation}
We may introduce a polynomial  $Q_{\alpha,l} \in \C[x]$ of degree $l$ such that
 \begin{equation} \label{eq:Laplacedas2}
 \mathcal{L}_{\R_+}\left(\zeta^{\alpha}\log(\zeta)^l\right)= \lambda^{-\alpha-1}Q_{\alpha,l}(\log(\lambda)).
 \end{equation}
Let $z\in \C$ with $\re(z)>0$ and let $m \in \N.$ We then have that 
\begin{align} \begin{split}
\label{Formal inverse}
\mathcal{L}_{\R_+}\circ \mathcal{B}(\lambda^{-z}\log(\lambda)^m)&=\lambda^{-z}\log(\lambda)^m, 
\\ \mathcal{B} \circ \mathcal{L}_{\R_+}(\zeta^{z-1}\log(\zeta)^m ) &=\zeta^{z-1}\log(\zeta)^m .
\end{split}
\end{align} 

\begin{lemma}  \label{lem:masterlemma}
Let $B: \R_+ \rightarrow \C$ be a measurable function and assume the integral defining $\mathcal{L}_{\R_+}(B)(\lambda)$ is absolutely convergent for $\re(\lambda)>0$. Assume there exists an increasing sequence $\{\alpha_j\}_{j=0}^{\infty}$ of real numbers strictly greater than $-1$ and a sequence $\{\beta_j \}_{j=0}^{\infty}$ of positive integers  giving an asymptotic expansion 
\begin{equation} \label{eq:tnearzero} B(t) \sim_{t \rightarrow 0} \sum_{j=0}^{\infty} b_{j} t^{\alpha_j} \log(t)^{\beta_j}.\end{equation}
Then the following holds  \begin{enumerate} \item There exists for large $ \lambda $ an asymptotic expansion $\tilde{\varphi}(\lambda)$ of the form
	\begin{equation}
	\mathcal{L}_{\R_+}(B) (\lambda) \sim_{\lambda \rightarrow \infty} \tilde{\varphi}(\lambda)
	\end{equation} 
	where
	\begin{equation} \label{eq:BorelA}
	\tilde{\varphi}(\lambda)  = \sum_{\alpha, \beta} b_{j} \lambda^{-\alpha_j-1}Q_{\alpha_j,\beta_j}(\log(\lambda))
	\end{equation} 
	and $Q_{\alpha_j,\beta_j} \in \C[x]$ is the degree $\beta_j$ polynomial introdced in \eqref{eq:Laplacedas2}.
	\item The Borel transform of $\tilde{\varphi}$ is equal to the expansion of $B$ 
	\begin{equation} \label{eq:BorelB}
	\mathcal{B}(\tilde{\varphi})(t)=  \sum_{j=0}^{\infty} b_{j} t^{\alpha_j} \log(t)^{\beta_j}.
	\end{equation}  \end{enumerate}
\end{lemma}

\begin{proof} The Lemma is an elementary consequence of Equations \eqref{eq:Laplacedas2} and  \eqref{Formal inverse}.   \end{proof}

The following theorem explains Borel-Laplace resummation. The content of Theorem \ref{thm:BorelLaplaceresummation} is standard in resurgence, and a proof can be found in e.g. \cite{Sauzin07}.

\begin{theorem} \label{thm:BorelLaplaceresummation} Let \begin{equation} \label{eq:DivergentSolution} \tilde{\varphi}(\lambda)= \sum_{j=0}^{\infty} c_j \lambda^{\alpha_j} \log(\lambda)^{\beta_j} \end{equation} be a formal series as in Definition \ref{def:Borel} with $\alpha_0>-1.$ Assume that \begin{itemize} \item  there exists a sector $S \subset \C$, such that for all $\zeta \in S$ of sufficiently small modulus the Borel transform $\mathcal{B}(\tilde{\varphi})(\zeta)$ converges to a holomorphic function $\hat{\varphi}(\zeta)$ (possibly upon choosing a  branch of $\log(\zeta)$ defined on $S$), and that
		\item  the function $\hat{\varphi}$ extends by analytic continuation along a half axis  $ \Gamma(\theta)=e^{2\pi i\theta}\R_+ \subset S$ (for some $\theta \in [0,1]$ say) and there exists a contant $C>0$ such that in a neighbourhood of $\Gamma(\theta)$ we have $\hat{\varphi}(z) = O(\exp(C \lvert z \rvert)).$
		
		\end{itemize}  Then the following holds.
	
	\begin{enumerate} \item  The Laplace transform $\mathcal{L}_{\Gamma(\theta)}(\hat{\varphi})$ is holomorphic on the open unbounded set $\{\lambda \in \C \mid \lvert \lambda \rvert >C,  \  \forall s \in S\setminus \{0\} \  \re(\lambda s)>0   \}.$
		\item The Laplace transform $\mathcal{L}_{\Gamma(\theta)}(\hat{\varphi})$  has $\tilde{\varphi}$ as its large $ \lambda $ asymptotic expansion
		\begin{equation}
\mathcal{L}_{\Gamma(\theta)}(\hat{\varphi})(\lambda) \sim_{\lambda \rightarrow \infty} \tilde{\varphi}(\lambda)
		\end{equation}
		\end{enumerate}
	\end{theorem}
One of the goal's of Ecalle's theory \cite{Ecalle81a,Ecalle81b} is to consider the case where the formal series $\tilde{\varphi}$  \eqref{eq:DivergentSolution} is obtained as a formal solution to some dynamical problem, which can be for instance an ODE or a difference equation (with a singularity at $\lambda^{-1}=0$). In such situations, the function $\mathcal{L}_{\Gamma(\theta)}\circ \mathcal{B}(\tilde{\varphi})$ will be a holomorphic solution, and resurgence is developed as a tool to analyze the monodromy (known as Stokes phenomena), which occur upon varying the choice of direction $\theta$ in which the Laplace transform is performed.

\end{appendices}

\bibliography{bibo}

\noindent 
       Jørgen Ellegaard Andersen\footnote{This work is supported in part by the center of excellence grant "Center for Quantum Geometry of Moduli Spaces" from the Danish National Research Foundation (DNRF95) and by the ERC-Synergy grant "ReNewQuantum".} \\
      Center for Quantum Mathematics\\
      Danish Institute for Advanced Study\\
        University of Southern Denmark\\
        DK-5000 Odense C, Denmark\\
     jea-qm{\@@}mci.sdu.dk
     \\\\
     William Elbæk Mistegård\footnote{This project has received funding from the European Union’s Horizon 2020 research and innovation programme under the Marie Skłodowska-Curie grant agreement No 754411.}  \\
      IST Austria \\
        AT-3400, Austria \\
    william.mistegaard{\@@}ist.ac.at

\end{document}